\documentclass{amsart}

\usepackage{amssymb}
\usepackage{url}
\usepackage[latin1]{inputenc}
\usepackage[T1]{fontenc}
\usepackage{graphicx}

\usepackage[labelfont=rm]{subfig}

\usepackage{array}
\usepackage{exscale,relsize}

\newcommand{\E}{\mathbb{E}}
\newcommand{\vv}{\mathbf{v}}

\usepackage{color}

\newcommand{\N}{\mathbb{N}}
\newcommand{\R}{{\mathbb R}}

\newcommand{\EE}{\mathbb{E}}
\newcommand{\X}{\mathcal{X}}

\newcommand{\Y}{\mathcal{Y}}
\newcommand{\Z}{\mathcal{Z}}
\newcommand{\M}{\mathcal{M}}
\newcommand{\one}{1\!\!1}

\newcommand{\bigoh}{\mathcal{O}}

\newcommand{\Ee}{{\mathrm{E}}}

\newcommand{\nklot}{S^{N-1}(\sqrt{N})}
\newcommand{\lpdist}{ {\mathrm d_{LP} }}

\newcommand{\empX}{\hat{\mu}_X}

\def\BB{{\mathcal B}}
\def\CC{{\mathcal C}}
\def\FF{{\mathcal F}}
\def\GG{{\mathcal G}}
\def\MM{{\mathcal M}}
\def\OO{{\mathcal O}}
\def\PP{{\mathcal P}}

\def\TT{{\mathcal T}}
\def\VV{{\mathcal V}}

\newcommand{\loi}{{\mathcal L}}

\def\SN{\mathfrak{S}_N}
\def\SSS{\mathfrak{S}}

\newcommand{\PPs}{{\mathcal P}_{\mbox{sym}}}

\newtheorem{thm}{Theorem}[section]
\newtheorem{prop}[thm]{Proposition}

\newtheorem{definition}[thm]{Definition}

\title[Random many-particle systems]{Random many-particle systems: applications from biology, and
  propagation of chaos in abstract models.}
\author{Bernt Wennberg}
\address{Department of Mathematics, Chalmers, SE41296 Gothenburg, Sweden \\
Department of Mathematics, University of Gothenburg, , SE41296 Gothenburg, Sweden}
\email{wennberg@chalmers.se}


\begin{document}

\maketitle

\begin{abstract}
 The paper discusses a family of Markov processes that represent many
 particle systems, and their limiting behaviour when the number of
 particles go to infinity. The first part concerns  model of
 biological systems: a model  for
 sympatric speciation, i.e. the process in which a genetically
 homogeneous population is split in two or more different species
 sharing the same habitat, and models for swarming animals. The second
 part of the paper deals with abstract many particle systems, and
 methods for rigorously deriving mean field models.

 These are notes from a series of lectures given at the 5$^{\mbox{th}}$
   Summer School on Methods and Models of Kinetic Theory, Porto
   Ercole, 2010. They are submitted for publication in "Rivista di
   Matematica della Università di Parma"  

\end{abstract}

  \section{Introduction}
  
  As the title suggests, these lecture notes consist of two rather
  different parts, although there is one uniting theme: random, interacting many
  particle systems. 
  
  The first part, dealing with applications from
  biology, begins with a model for {\em sympatric speciation}; this is
  the process in which a population of animals or plants is split in two
  or more separate species, remaining in the same geographical area
  (hence the word {\em sympatric}). In this case, the particles are
  individuals, and the interaction is the mating procedure and the
  selection process. This part is based mostly on ~\cite{HenrikssonLundhWennberg2010}.
  A rather different class of models, more similar to
  the classical kinetic theory of gases, are models for swarms (that
  could be swarms of insects, flocking birds, schooling fish, or for
  that matter, crowds of people). The particles are then the
  individuals, and the interaction is usually the voluntary motion of
  the individuals, based on their visual perception of other individuals
  in the neighborhood. Such problems have attracted a lot of interest
  in the kinetic theory community recently, and although I will give some
  of the important references, these notes do not  give a
  complete review of the current works, only to give another example of
  how ideas from kinetic theory can be applied to biological
  problems. To a large extent it is based on ongoing research with Eric
  Carlen and Pierre Degond~\cite{CarlenDegondWennberg2011}.
  
  The remaining part of the notes deal with {\em propagation of chaos},
  which very vaguely means that if the particles initially are distributed
  independently of each other in phase space, then they remain
  independent along the evolution of the system. This never holds for
  interacting particle systems of the kind considered here, as long as
  the number of particles is finite, but for some models it can be
  proven to hold in the limit of infinitely many particles. It is one of
  the major challenges in kinetic theory to prove that propagation of
  chaos holds for real particle systems, a topic that is discussed in
  detail in Pulvirenti's notes in this
  issue~\cite{PulvirentiPortoErcole2010}. Here we discuss a
  much easier case, where the microscopic model is already random, a
  family of Markov jump processes in state spaces $\E^n$ that represent
  $n$-particle configurations, and the corresponding master
  equations. The propagation of chaos can then be expressed in terms of
  the marginals of the $n$-particle distributions.  This approach to the
  propagation of chaos 
  goes back to Mark Kac~(\cite{Kac1957}), and the key ideas in that paper
  will be presented. A different approach was taken by
  Grünbaum~(\cite{Grunbaum}), closely related to de Finetti's
  theorem on the conditional independence of exchangeable
  observations. This approach has been taken a step further
  in~\cite{MiMoWe2011arXiv}, from which much of the material in these notes
  is taken. And~\cite{MiMoWe2011arXiv} was inspired in part by the lectures
  of P.L.Lions on mean field games~(\cite{PLL_cours_champsmoyens}). A
  small section of these notes is essentially taken from one of the 
  first lectures in his series.
  
  What unites these rather disperse topics is that they all deal with 
  Markov processes in a state space $\E^n$, {\em i.e.} an $n$-fold
  product of an Euclidean space $\E$, or some  submanifold
  of $\E^n$ representing {\em e.g.} the conservation of energy.\footnote{It is not necessary that
    $\E$ be Euclidean, it should be a Polish space -- a separable,
    completely metrizable topological space} In the state
  $\vv=(v_1,...,v_n)\in \E^n$, each component $v_k$ represents the state
  of one particle. In Kac's original model, and other models that
  represent real gases, the jumps only change two components, $v_i$ and
  $v_j$, say, simultaneously, although the rate at which a particular
  couple of particles interact may be determined as a function of the
  full state. No deep knowledge of Markov processes is needed to read
  these notes, only a basic understanding of the definitions is
  assumed. A comprehensive book on the topic is~\cite{EthierKurtz}, and a standard
  reference with applications in physics and chemistry
  is~\cite{vanKampen2007}.

  \section{Applications from biology: a model for sympatric speciation}
  
  There are several mathematical models for speciation that are related
  to the models from the kinetic theory of gases. The one that is
  presented here comes from~\cite{HenrikssonLundhWennberg2010}, but
  there are many other examples, and I will very briefly mention a
  couple. But first we need to reflect over the concept of
  species. Although most of us have a vague idea of what a species is,
  it is by no means an easy task to make a proper definition. Until at
  least the 18th century, the flora and fauna were thought of as being
  rather stationary, and a species was characterized by producing an
  offspring of (essentially) the same kind. Notably Linnaeus created a
  taxonomic system for classifying and naming the species, a system that
  is still used today. But it does not really define the concept of a
  species, rather it 
  gives a hierarchical structure, with similar plants, or animals, grouped
  together. Compte de Buffon, contemporary with Linnaeus
  (both of them were born in 1707), classified two individuals as belonging
  to the same species if they can produce fertile offspring. This
  definition is problematic for several reasons, one being that it is
  not a transitive relation: One could have three candidates for a
  species, A, B, and C, such that A and B can produce fertile offspring,
  B and C too, but not C and A. It may also happen that the result
  depends on whether A or B is female. The discovery of DNA and
  techniques to analyze the genetic code has provided new means for
  classifying species, but there is no general definition of ``species''
  that is useful in all situations.
  
  With Darwin's {\em On the origin of species}~\cite{Darwin}, a mechanism for
  evolution was described: due to phenotypic variation within a
  population, some individuals will reproduce less efficiently, and
  there will be a selection against this phenotypic character. But this
  mechanism is not enough to explain how a species can evolve into two
  different species. A nice discussion on this topic can be found in the
  introduction to van Doorn's thesis {\em Sexual selection and sympatric
    speciation\,}\cite{vanDoorn_thesis}. 
  
  {\em Allopatric speciation}  may happen if a homogeneous population is split
  into two geographically separated regions, such as two islands. By
  selection the two sub-populations will then evolve to adapt to the
  local environment, but also phenotypic characters that are not
  selected against will also change, and eventually the two
  sub-populations may be so different that they have become two different
  species. It is much more difficult to understand {\em sympatric
    speciation}, where the two sub-population share the same geographical
  area. van Doorn lists a number of obstacles for speciation to take
  place, and exemplifies this with birds feeding on different size
  grains: small, medium and large. The fitness of a bird is quantified
  in terms of the feeding rate. Speciation would now mean, for example,
  that one sub-population specializes in feeding on small seeds and
  another one on large seeds, but for this to happen, the population
  must reach a state known as ``disruptive selection'', {\em i.e.} a
  situation where the feeding rate could improve by
  changing a phenotypic character (such as the beak length)  {\em
    either } in one  direction {\em or} the other. A first step
  towards speciation is taken 
  if two sub-populations evolve in different directions, leading to
  a phenotypic ``polymorphism'', but unless the ecological landscape
  gives an advantage to the smaller sub-population, only the larger one
  will remain, and hence the polymorphism is eventually lost, and no
  speciation can take place. 
  
  The next step towards speciation is the evolution of a ``reproductive
  isolation'', a mechanism that prevents the formation of hybrids. A key
  concept is that of ``assortative mating'', meaning that reproduction
  takes place essentially within sub-populations: individuals choose
  mating partners according to specific criteria. Finally, some kind of
  dependence should develop between the genes that are responsible for
  the fitness to the landscape (beak length, in our example), and the
  genes responsible for the assortative mating.
  
  All this is far from completely understood, and the literature on
  the subject is   vast. The model for sympatric speciation presented
  here is one example that addresses   the question.
  
  \subsection{Adaptive dynamics}
  
  One approach to evolution and speciation is {\em adaptive dynamics}. A
  recent book treating this is~\cite{DercoleRinaldi2008}, which states
  that adaptive dynamics is ``the long term evolutionary dynamics of
  quantitative characters driven by the processes of mutation and
  selection'', and is a theory that has been developed for example by
  Geritz et al. \cite{Geritz_etal1998}. A short introduction is given
  in~\cite{hitchhikerAD}.
  
  In the easiest case, adaptive dynamics is concerned with the evolution
  of a scalar trait in a monomorphic population. A scalar trait could be
  for example the length of a bird's beak, $x>0$, say, and that the
  population is 
  monomorphic means that all individuals have exactly the same beak
  length. Adaptive dynamics takes place on a time scale much longer than
  the typical time scale of a population, and hence it is assumed that
  the resident population with beak length $x$ is stationary. The main
  issue of adaptive dynamics is to understand what happens to a small
  group of individuals that (by mutation) has a different trait value,
  $y$, say. The initial growth rate of the rare mutant population,
  denoted $S_x(y)$, is sometimes called the invasion exponent. Because
  the resident population is assumed to be stationary, $S_x(x)=0$. If
  $S_x(y)>0$, then the mutant is more fit, and will eventually take
  over, whereas if $S_x(y)<0$, the mutant will not invade. The {\em
    selection gradient} determines the direction of evolution of the
  trait: if  $S_x'(x)>0$, then an invading mutant with trait $y>x$ has a
  better fitness, and will replace the resident population. The new
  resident population will have trait value $y$. Similarly, if
  $S_x'(x)<0$, then the resident population will be replaced by a
  population with smaller trait value. Values of $x$ such that
  $S_x'(x)=0$ are known as {\em evolutionarily singular states}, and if
  it corresponds to a local fitness maximum, it is called an {\em
    evolutionary stable strategy}. A resident population at an ESS
  cannot be invaded by a nearby mutant, because all nearby strategies
  are less fit to the environment. {\em Disruptive selection} can occur
  when an evolutionary singular state is at a local fitness minimum. The
  a nearby mutation at either side of $x$ has a better fitness, and this
  is what is needed for speciation to take place.
  
  Another concept is a {\em convergence stable strategy}, which is a
  strategy such that monomorphic populations with $y$ close to $x$ can
  be invaded by mutants which are even closer to $x$. When this is the
  case, the trait value of the resident population will approach $x$.
  
  The mutations arrive in a population randomly, and it is not
  necessarily true that the mutants have trait values close to that of
  the resident population, but if the mutants are small, and the
  frequency of mutations is scaled properly, it is possible to derive an
  ODE which describes the rate of change of the trait of the resident
  population, the {\em canonical equation of adaptive dynamics}.
  
  An approach based on the Hamilton - Jacobi equations can be found
  in~\cite{DiekmannJabinMischlerPerthame2005}. 
  
  \subsection{Examples of mathematical models of speciation}
  
  There are several examples of mathematical models for the competition
  within a population structured according to some phenotypic
  trait. Desvillettes {\em et
    al.}~\cite{DesvillettesJabinMischlerRaoul2008}  consider the
  following model of logistic type:
   \begin{equation*}
     \frac{\partial f}{\partial t}= \left(a(y)-\int_{Y} b(y,y')
       f(y')dy' \right)f\,. 
   \end{equation*}
   Here $f(y)$ is a density describing the distribution of the
   population according to the trait $y\in\Y$, where $\Y$ is
   compact. The birth rate of individuals with trait $y$ is $a(y)$, and
   death rate is 
   \begin{equation*}
     \int_{\Y} b(y,y') f(y')dy'\,.
   \end{equation*}
  The death rate can be seen as a model for competition within the
  population, the function $b(y,y')$ giving death rate of an individual
  of trait $y$ due to the interaction with an individual of trait $y'$. 
  
  The authors prove global in time existence and uniqueness in
  $L^1(\Y)$ for this equation, assuming sufficient regularity of
  the functions $a(y)$ and $b(y,y')$. They also present the results of
  numerical simulations that show how an initially unimodal trait
  distribution evolves into a bimodal, and then multimodal
  distribution. In fact, they even show that a limiting solution must
  consist of a sum of Dirac masses.
  
  A different model is provided by Méléard and
  Tran~\cite{MeleardTran2009}, where an age structured population is
  studied (this paper is an extension of earlier works by Méléard and
  co-authors, see the reference list of~\cite{MeleardTran2009}). In their
  model the population is described by a random measure,
  \begin{equation*}
    Z_t(dx,da) = \sum_{i=1}^{\langle Z_t,1\rangle} \delta_{(x_i(t),a_i(t))}\,.
  \end{equation*}
  The size of the population is $\langle Z_t,1\rangle $, and each
  individual is characterized by its trait value $x\in\X$ and its age $a\in\R^+$.
  Each individual produces offspring with rate $b(x,a)$ depending on the
  trait $x$ and age $a$, and the offspring is born with age $a=0$ and a
  trait $x=x_{par} + h$, {\em i.e.} the parent's trait plus a mutation
  difference $h$ which is random and distributed with law
  $k(x,h)dh$.\footnote{A natural variation of this would be to consider a
    mutation rate also depending on the parent's age} Just like
  in~\cite{DesvillettesJabinMischlerRaoul2008} the death rate has a
  component due to competition between all individuals in the
  population:
  \begin{equation*}
  d_{tot}=  d(x,a) + \int_{\R^+\times\X}   U((x,a),(y, \alpha)) Z(dy,d\alpha)
  \end{equation*}
  In simulations, using birth and death rates
  \begin{eqnarray*}
   b(x,a) &=& x( 4-x)e^{-a}\qquad\mbox{and}\\
   d_{tot}(x,a)&=  & \frac{1}{4} + a
       \int_{\tilde{\mathcal{X}}} 
       C\left(1-\frac{1}{1+\nu \exp(-k(x-y))} \right) Z(dy,d\alpha)\,,
  \end{eqnarray*}
  they find that an initially monomorphic population may evolve into a
  population with a bi-modal trait distribution. However, one of the main
  objectives of their paper is to study the ``large population -- rare
  mutation''-scaling. Setting
  \begin{eqnarray*}
     Z^n_t = \frac{1}{n}\sum_{i=1}^{n\langle Z^n_t,1\rangle}
       \delta_{(x_i(t),a_i(t))}\,,
  \end{eqnarray*}
  they prove that  $Z^n_t \rightarrow \xi_t \in \mathcal{M}_F$, where
  $\mathcal{M}_F$ is the set of finite measures on
  $\R^+\times\X$. Actually,    
  $ \xi \in \mathcal{C}(\R_+,\mathcal{M}_F) $
  and, for all  $ t>0,
  f\in\mathcal{C}^{0,1}_b(\tilde{\mathcal{X}},\R)$
  \begin{eqnarray*}
    \lefteqn{
      \langle\xi_t,f\rangle = \langle \xi_0,f\rangle  +}\nonumber \\
    & &  \int_0^t \int_{\tilde{\mathcal{X}}}\left[ \partial_a
      f(x,a)+f(x,0)b(x,a)-
      f(x,a)(d(x,a)+ \xi_s U(x,a))\right] \xi_s(dx,da)\,.     
  \end{eqnarray*}
  
  A model that in many ways is similar to the one that is presented in
  the next section can be found in a paper by Dieckmann and
  Doebeli~\cite{DieckmannDoebeli1999}. Actually they discuss two
  different models, of which one is an individual based simulation
  model, the other a model in the framework of adaptive dynamics. The
  resident population,  having phenotype $x$,  is assumed to satisfy the
  following logistic equation:
  \begin{equation*}
    \frac{dN(x,t)}{dt} = r N(x,t) \left[
         1-\frac{N(x,t)}{K(x)}\right] \,,
  \end{equation*}
  where $N(x,t)$ is the size of the population at time $t$, and 
  $K(x)$ is the carrying capacity for a monomorph population with trait
  $x$. In~\cite{DieckmannDoebeli1999} $K(x)$ is a Gaussian with mean
  $x_0$ and variance $\sigma_K^2$. Due to competition with the resident
  population, a rare mutant with trait $y$ will grow with rate
       $\displaystyle r \left[1-\frac{C(x-y)K(x)}{K(y)}\right]$; here
       $C(x-y)$, which describes the strength of the competition between
       a phenotype $x$ and a phenotype $y$, is a Gaussian with variance
       $\sigma_C^2$. The mathematical analysis
       in~\cite{DieckmannDoebeli1999} shows that evolutionary branching
       can take place only if $\sigma_C<\sigma_K$.

  \subsection{A model of sympatric speciation through
  reinforcement}
  
  As explained in the introduction, even if evolution brings the
  resident population 
  to a state of disruptive selection, it is often more natural for the
  population to evolve in one direction rather than to split in two
  viable sub-population evolving in different directions. If the latter
  is to happen there must be some mechanism to favor the smaller of the
  sub-populations. Competition of resources, and specialization to
  particular parts of the available resources may be one such
  mechanism. For the sub-populations to evolve into two different
  species, some kind of reproductive isolation is needed to prevent the
  formation of hybrids. In~\cite{HenrikssonLundhWennberg2010} we have
  developed a model to study some aspects of this.
  
  We recall that {\em sympatric speciation} means that a population
  develops into two species sharing the same geographical habitat (but
  usually not sharing the same resources). By {\em reinforcement} we
  mean a process by which natural selection strengthens the separation
  of the sub-populations. In this model, reinforcement is implemented
  via a characteristic trait $y$ describing the appearance of an
  individual ({\em e.g.} the color of the tail feathers, the pitch of
  the song ...), and another trait $y^*$ describing what characteristic
  traits in a potential mating partner the individual is attracted
  to. We assume that the traits $y$ and $y^*$ are only related to the
  choice of partner, and not directly to the fitness. Fitness, on the
  other hand, is determined by a trait $x$, which is related to the
  distribution of food resources in the local environment. We show, by
  simulation, that in this model reinforcement is needed for speciation
  to take place, and that the expected time before the speciation event
  is shorter if the characteristic trait $y$ has more than one dimension. 
  
  Here is the model:
  
  The population lives in an environment, where food (or other essential
  resources) is characterized by a parameter $x\in\X$, and that the food
  is distributed in space according to a density $f(x)$. The trait of an
  individual in the population has several parts:
  \begin{itemize}
  \item $x\in X$ is related to the fitness, its competitivity in
    collecting the essential resource;
  \item $y\in\Y$ is a recognizable, characteristic trait, and $y^*\in\Y$
    is the individual's preference of trait value in potential mating
    partners. These two traits combined yield the probability that a
    given couple of individuals will mate.
  \end{itemize}
  
\begin{figure}
    \centering
    \includegraphics{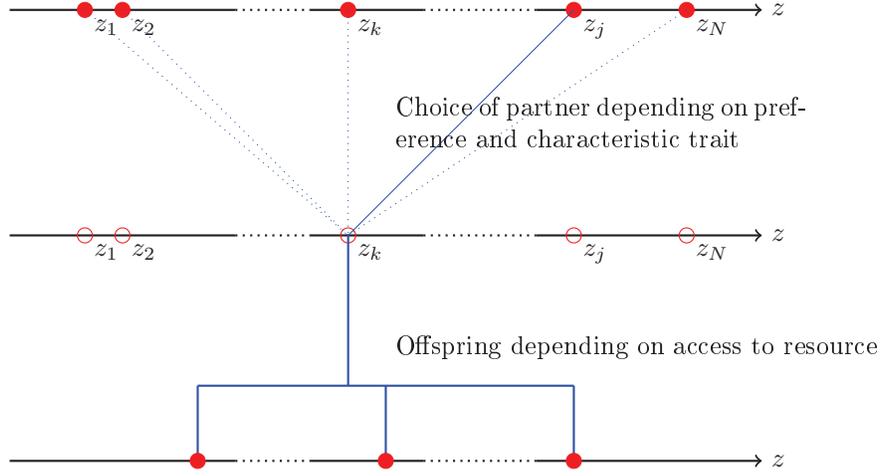}
    \caption{A pictorial description of the replicator dynamics of the population}
    \label{fig:1}
  \end{figure}

  The size of the population is denoted $N$, and letting $z_k
  =(x_k,y_k,y^*_k)\in\Z=\X\times\Y\times\Y$ be the phenotype of individual $i$,   
  the whole phenotype distribution in the population is
  \begin{equation*}
    Z^N = (z_1,...,z_N)\in\Z^N\,.
  \end{equation*}
  The dynamics is time-discrete, and we assume that only the offspring
  survives from one generation to the next. The process can be described
  as follows:
  \begin{enumerate}
  \item Each individual collects food according to its relative fitness,
  \item It then chooses a mating partner, at random but with a high
    probability to select a partner with a characteristic trait
    corresponding to the preference.
  \item The size of the offspring is Poisson distributed with a
    parameter proportional to the couple's joint access to the food
    resource.
  \item The phenotype of the offspring is the average of that of the
    parents', but mutations are included by adding a (Gaussian) random variable.
  \end{enumerate}
  The procedure is described in Figure~\ref{fig:1}.
  More precisely
  \begin{enumerate}
  \item An individual $k$ has access to  a fraction $\ c_k $ of the
    available resource:
    \begin{eqnarray*}\
    c_k= \int_{\X} 
      \frac{ e^{-(x_k-x)^2/2\sigma_x^{2}} }
        {\sum_{j=1}^{N_t}  \int_{\X} e^{-(x_j-x)^2/2 \sigma_x^2}}\,.
    \frac{df(x)}{\|f\|}
  \end{eqnarray*}
  This represents the {\em competition} among the individuals.
  \item Each individual $k$ is given the opportunity to choose a
    mating partner, and 
    chooses $j$ with probabilily
   \begin{eqnarray*}
      \mbox{Prob}(j_k=j) &= &\left\{  
        \begin{array}{ll}
          \frac{ e^{-|y^*_k-y_j|^2/2\sigma^2}}{ \sum_{i \ne k}
            e^{-|y^*_k-y_i|^2/2\sigma^2}}&\qquad(j \ne k)\,,\\
          \\
          0 &\qquad (k=j)\,.
        \end{array}\right.
      \end{eqnarray*}
  This is the {\em reinforcement} in the model, because it helps forming
  sub-populations such that mating takes place within the group. The
  parameter $\sigma$, which we have taken to be  the same for all
  individuals in the population, determines the choosiness in selection
  of partners for mating. 
  \item The couple $(k,j_k)$ then produces a Poisson distributed number
    $n_{\kappa}$ of children, with rate $\frac{c_k+c_{j_k}}{2}\|f\|$,
    {\em i.e.} proportional to the amount the the resource that has
    been collected by the couple. This means that the size of the
    population at time $t+1$ will be a Poisson distributed variable,
    $N_{t+1}=\sum_{k=1} n_k$ with a random parameter
    \begin{equation*}
      \lambda = \|f\| \left(1+\sum_{k=1}^{_t} c_{j_k}\right).
    \end{equation*}
   It is random because of the random choice of partner, $j_k$, and the
   law depends on the whole population at time $t$.
  \item  Each child has a trait $ z=z_{k,i}$, $ i=1...n_k$
    \begin{eqnarray*}
      z_{k,i} &=& \frac{z_k+z_{j_k}}{2} + (\xi_i,\eta_i,\eta^*_i)\,.
    \end{eqnarray*}
  
      where $(\xi_i,\eta_i,\eta^*_i)$ are Gaussian random variables.
  \end{enumerate}

\begin{figure}[h!]
    \centering
  \subfloat[]{\label{fig:runA_BW}%
    \includegraphics[width=0.38\textwidth]{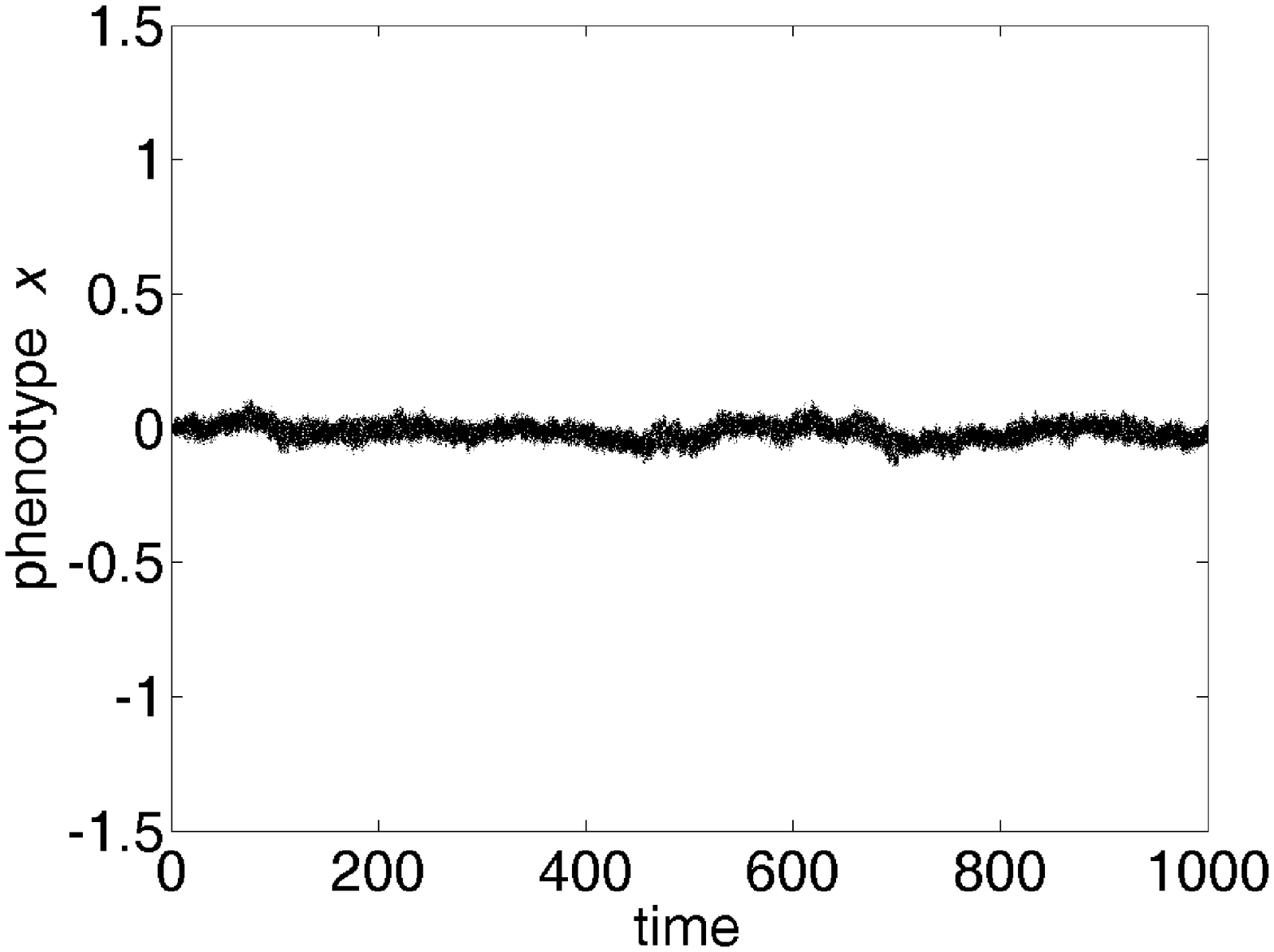}}\hspace{0.1\textwidth}
  \subfloat[]{\label{fig:runB_BWa}%
    \includegraphics[width=0.38\textwidth]{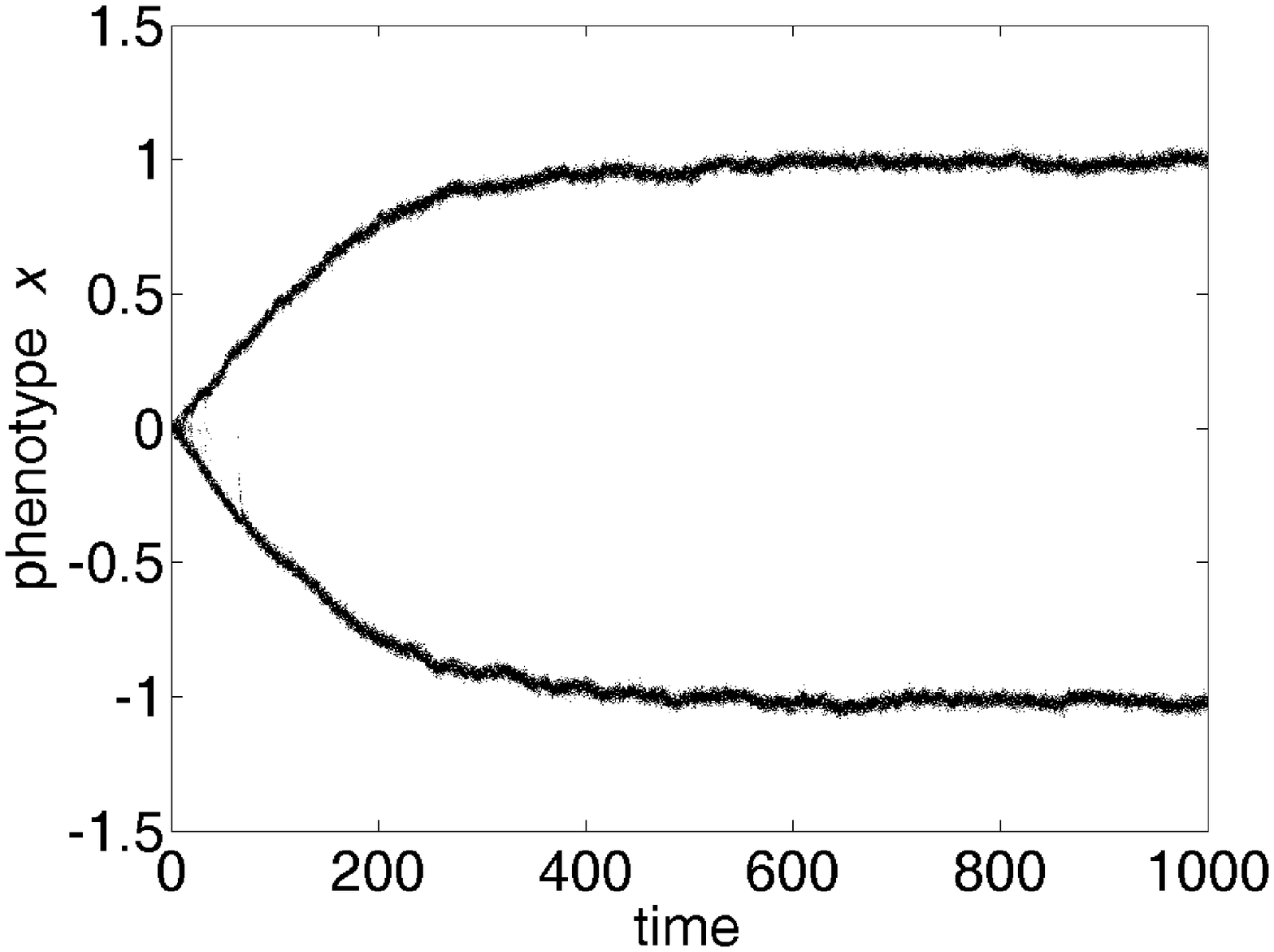}}
  
  \subfloat[]{\label{fig:runB_BWc}%
    \includegraphics[width=0.38\textwidth]{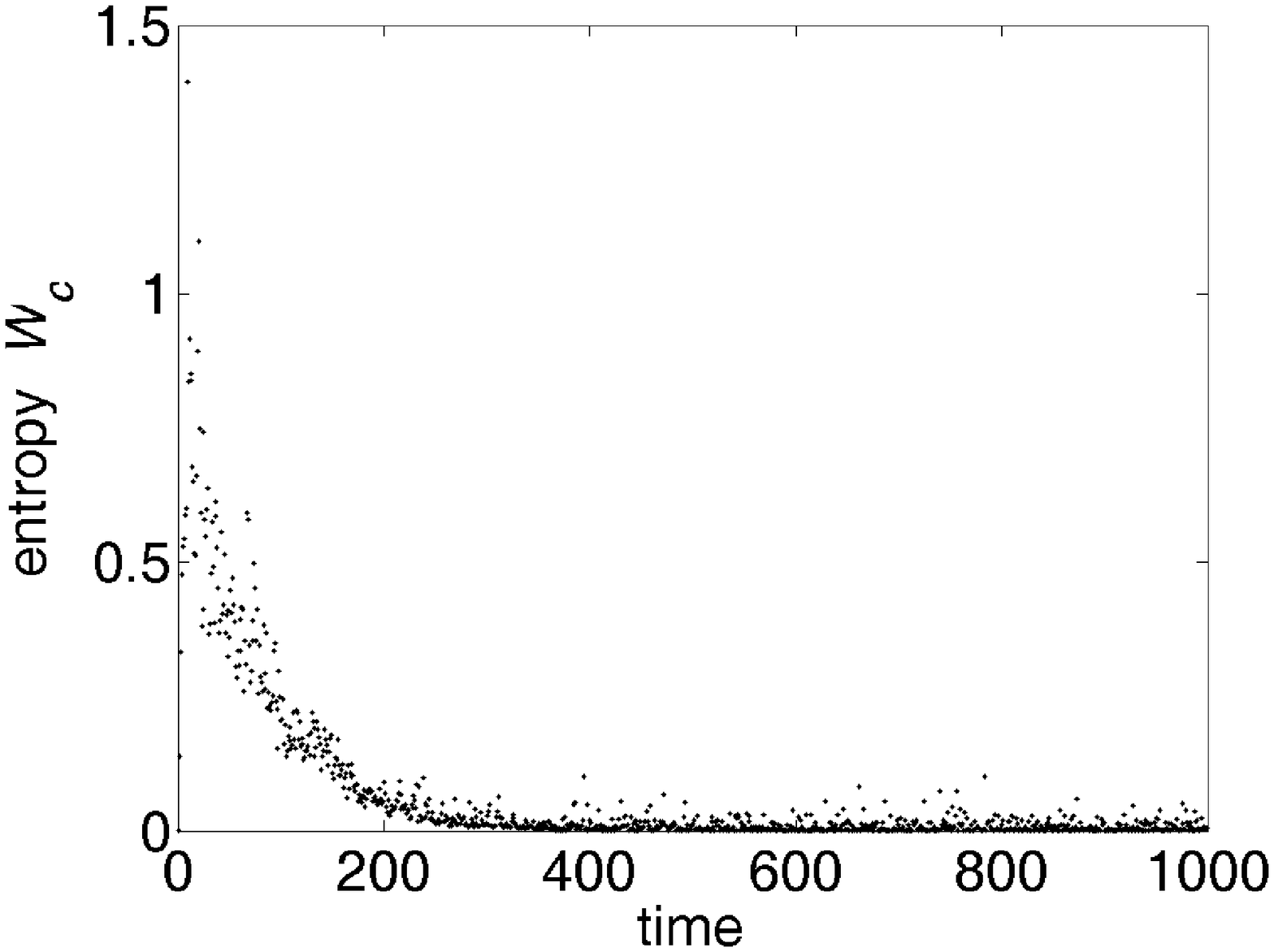}}\hspace{0.1\textwidth}
  \subfloat[]{\label{fig:runB_BWb}%
    \includegraphics[width=0.38\textwidth]{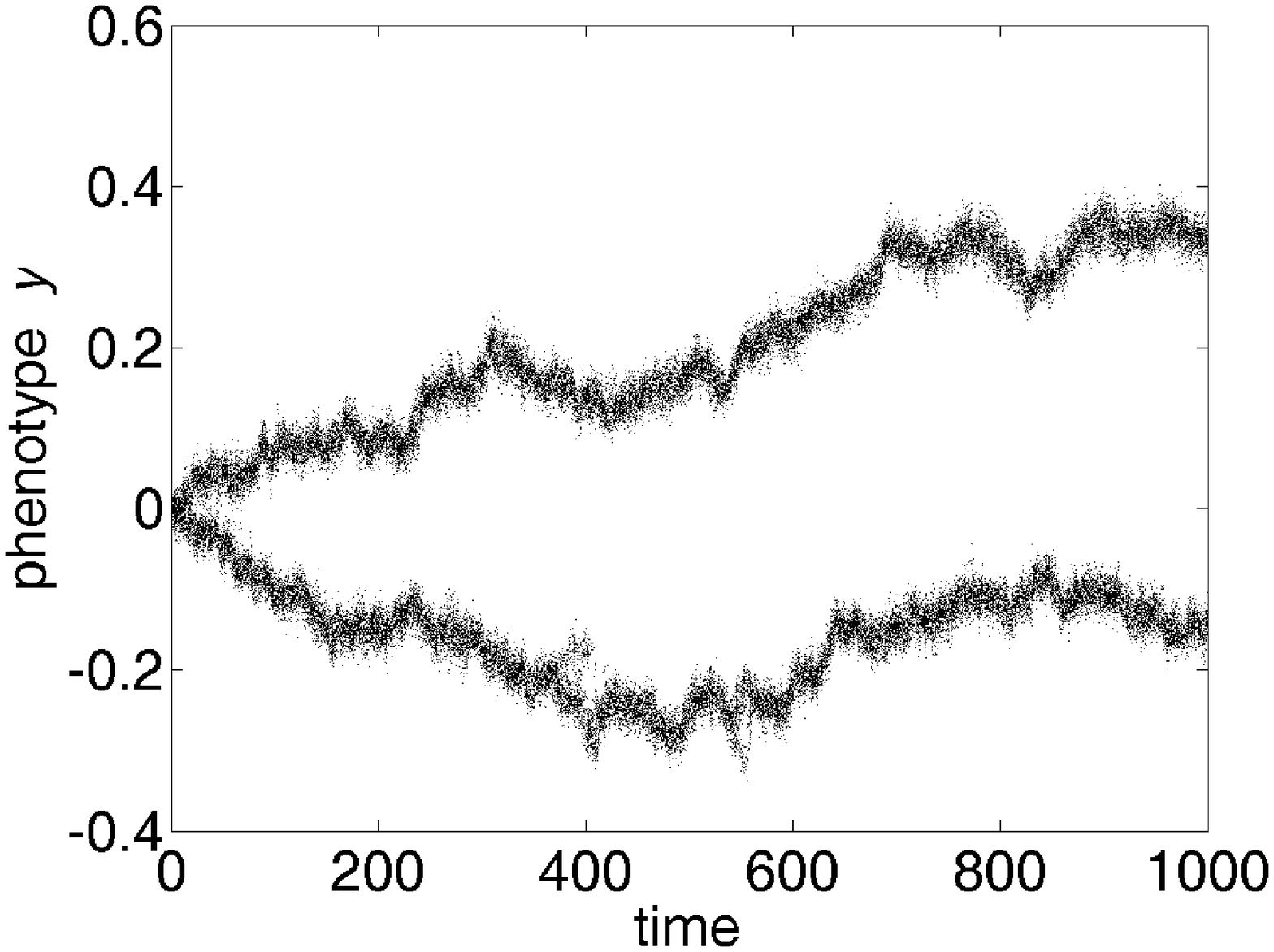}}
  
    \caption{(a) and (b)  $x$-trait without and with reinforcement. \\
           (c) $c$-entropy, $W_c$, (d):  appearance trait $y$}
  \label{fig:runAB}
  \end{figure}

\begin{figure}[h!]
    \centering
  \subfloat[]{\label{fig:runC_BWa}%
    \includegraphics[width=0.38\textwidth]{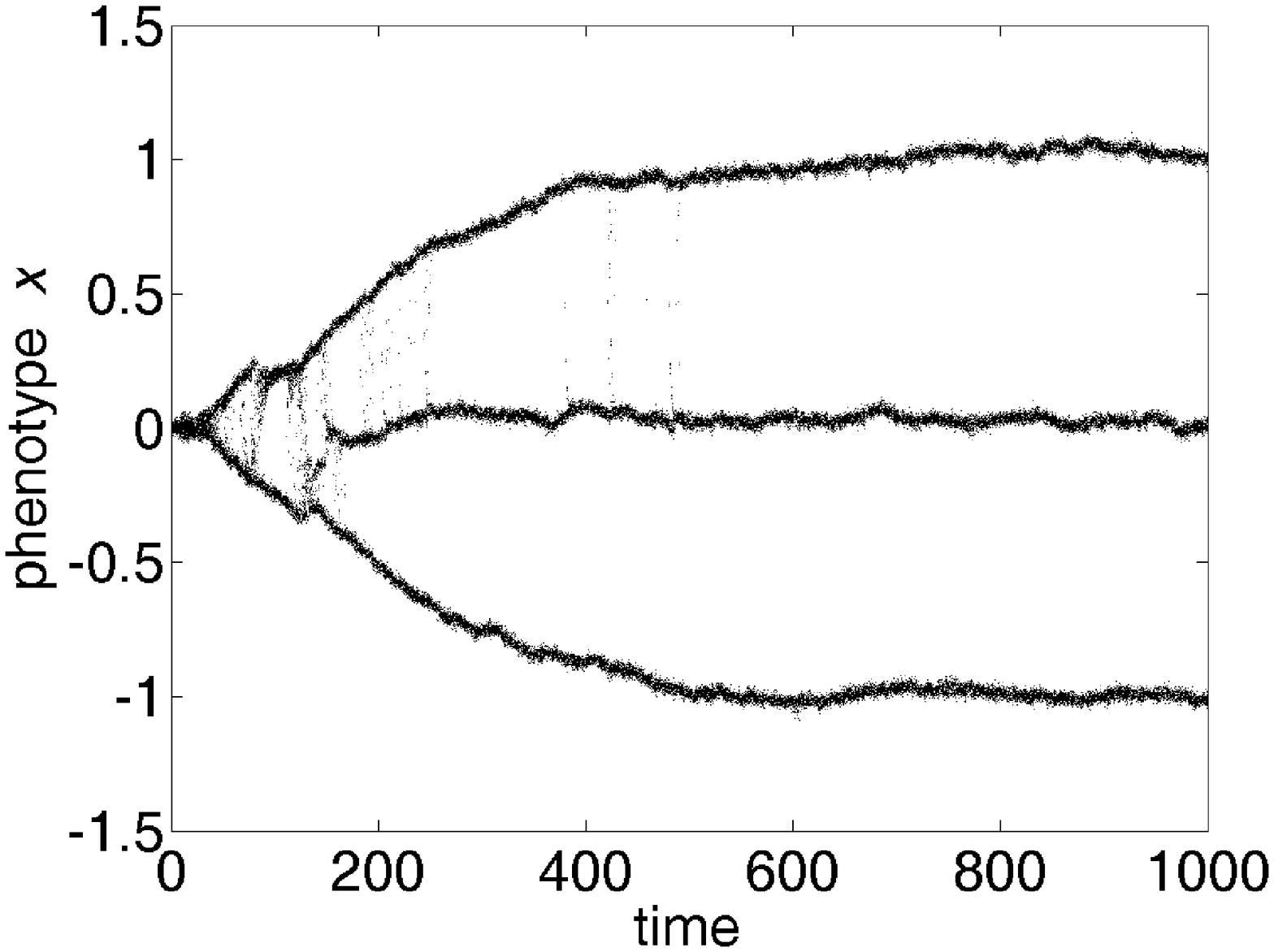}}\hspace{0.1\textwidth}
  \subfloat[]{\label{fig:runC_BWb}%
    \includegraphics[width=0.38\textwidth]{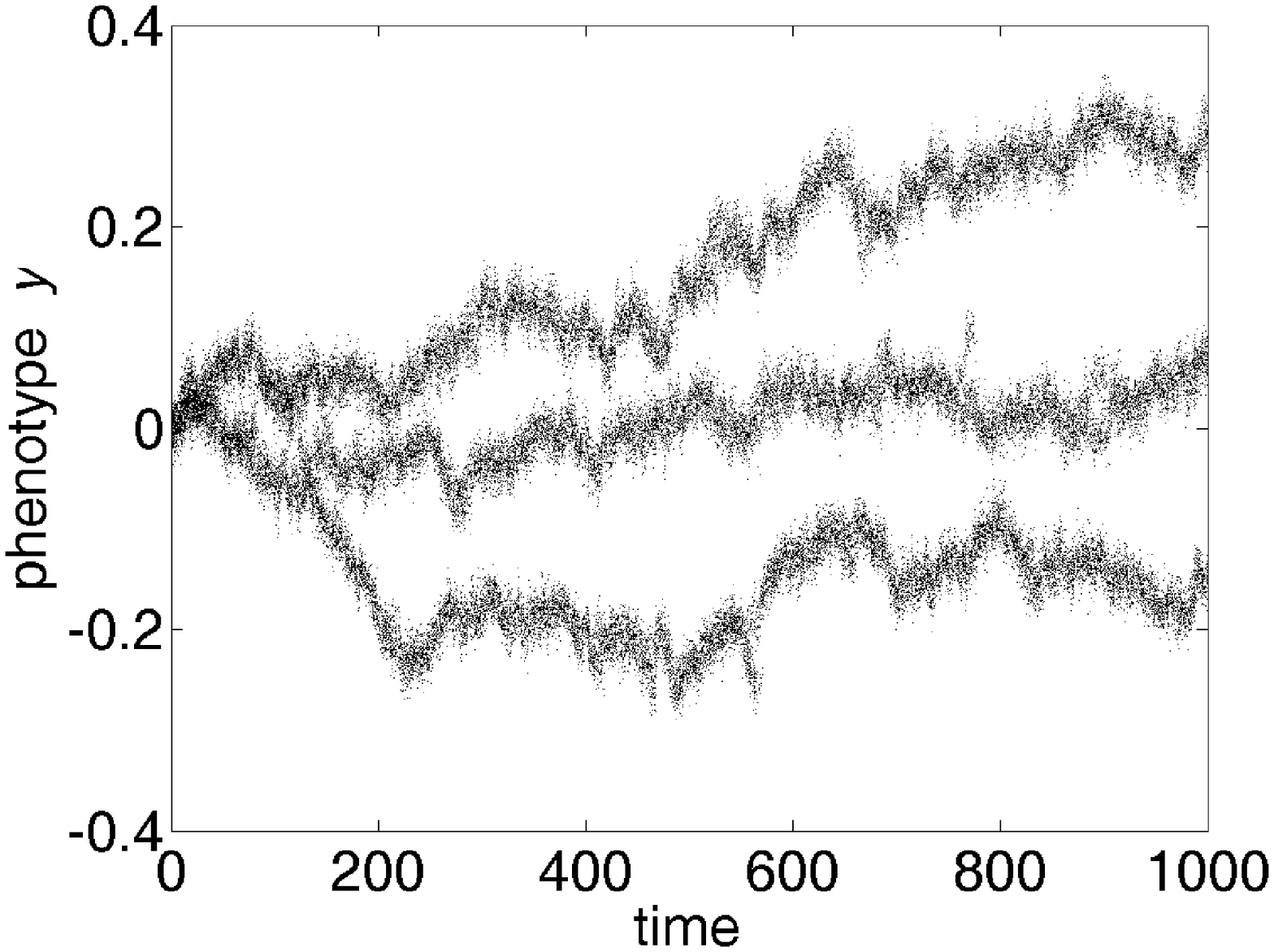}}
  
  \subfloat[]{\label{fig:runC_BWc}%
    \includegraphics[width=0.38\textwidth]{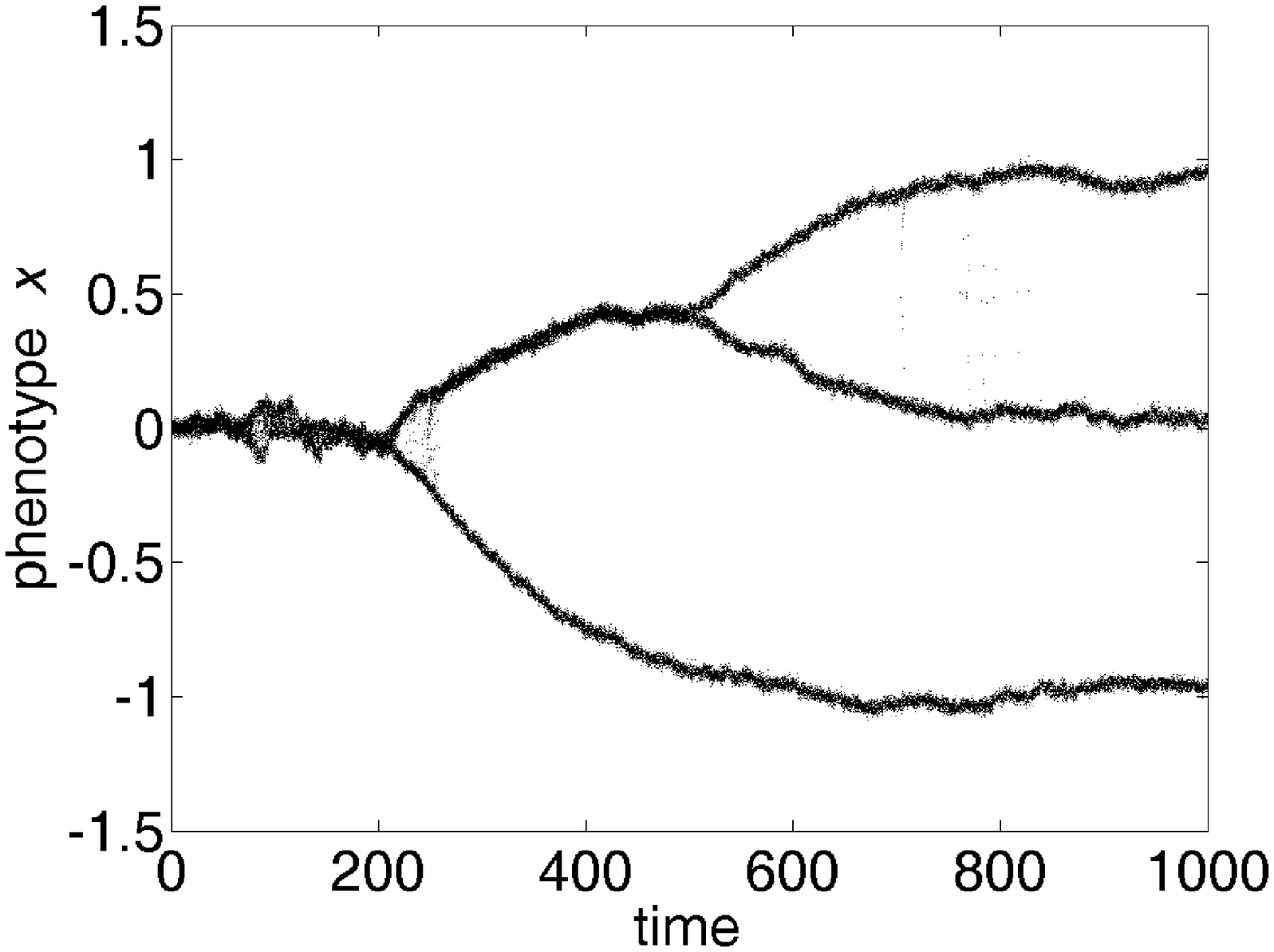}}\hspace{0.1\textwidth}
  \subfloat[]{\label{fig:runC_BWd}%
    \includegraphics[width=0.38\textwidth]{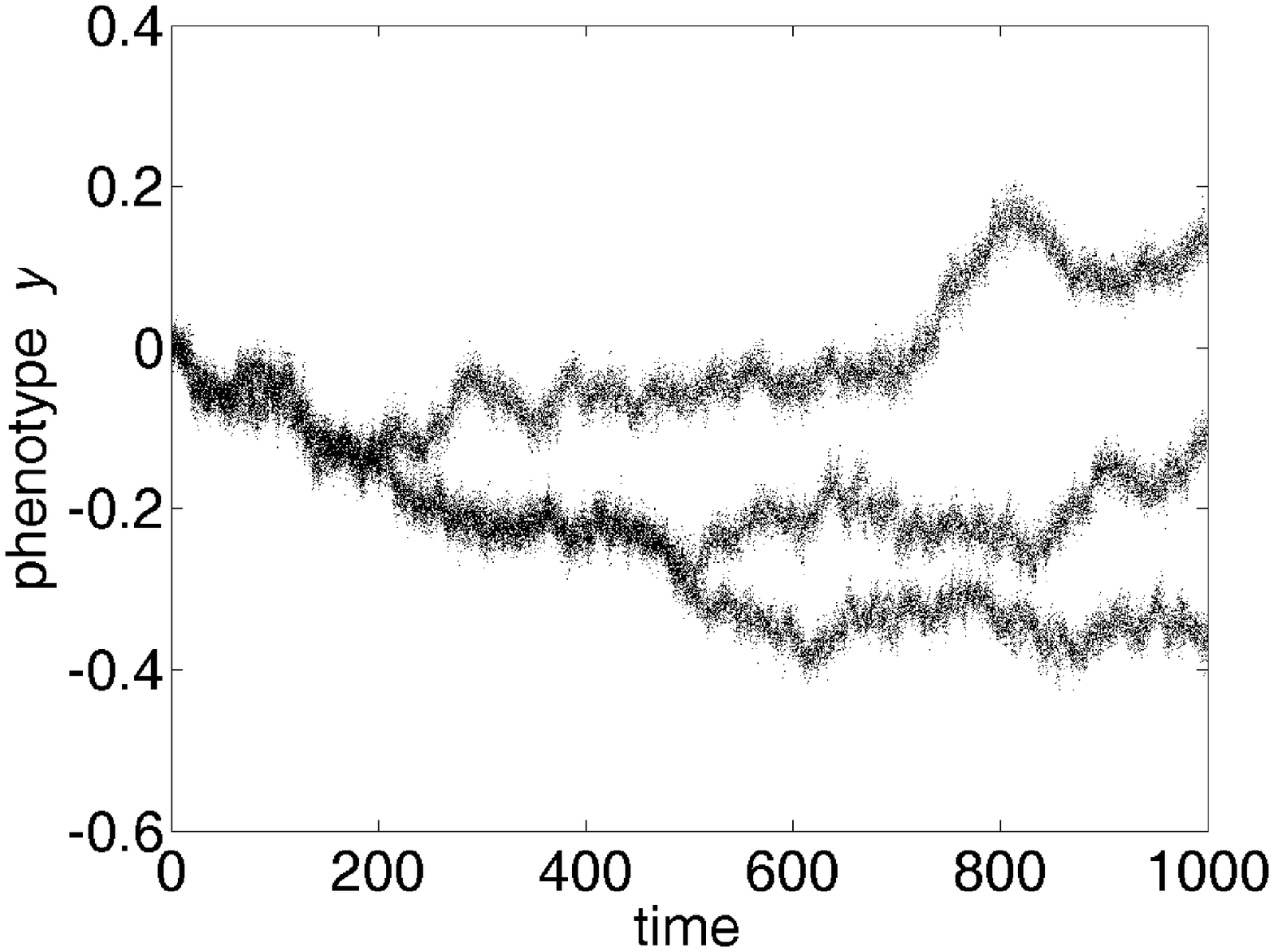}}
  
    \caption{Two simulations with identical  parameters. \\
      (a) and~(c): trait $x$ in the population, and~(b) and~(d):  trait $y$. }
  \label{fig:runC}
  \end{figure}

\begin{figure}[h!]
    \centering
  \subfloat[]{\label{fig:sim_x3a}%
    \includegraphics[width=0.38\textwidth]{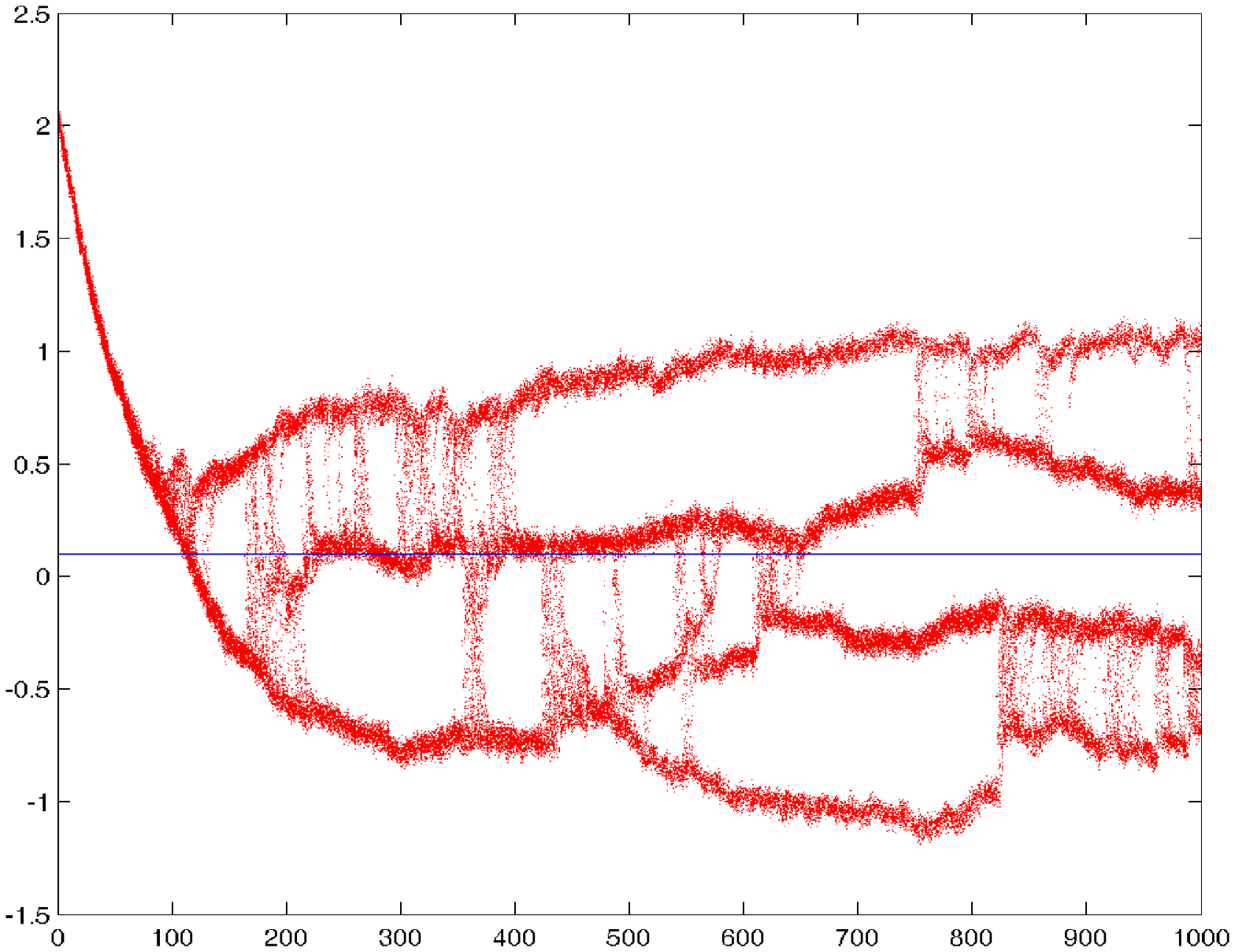}}\hspace{0.1\textwidth}
  \subfloat[]{\label{fig:sim_x3b}%
    \includegraphics[width=0.38\textwidth]{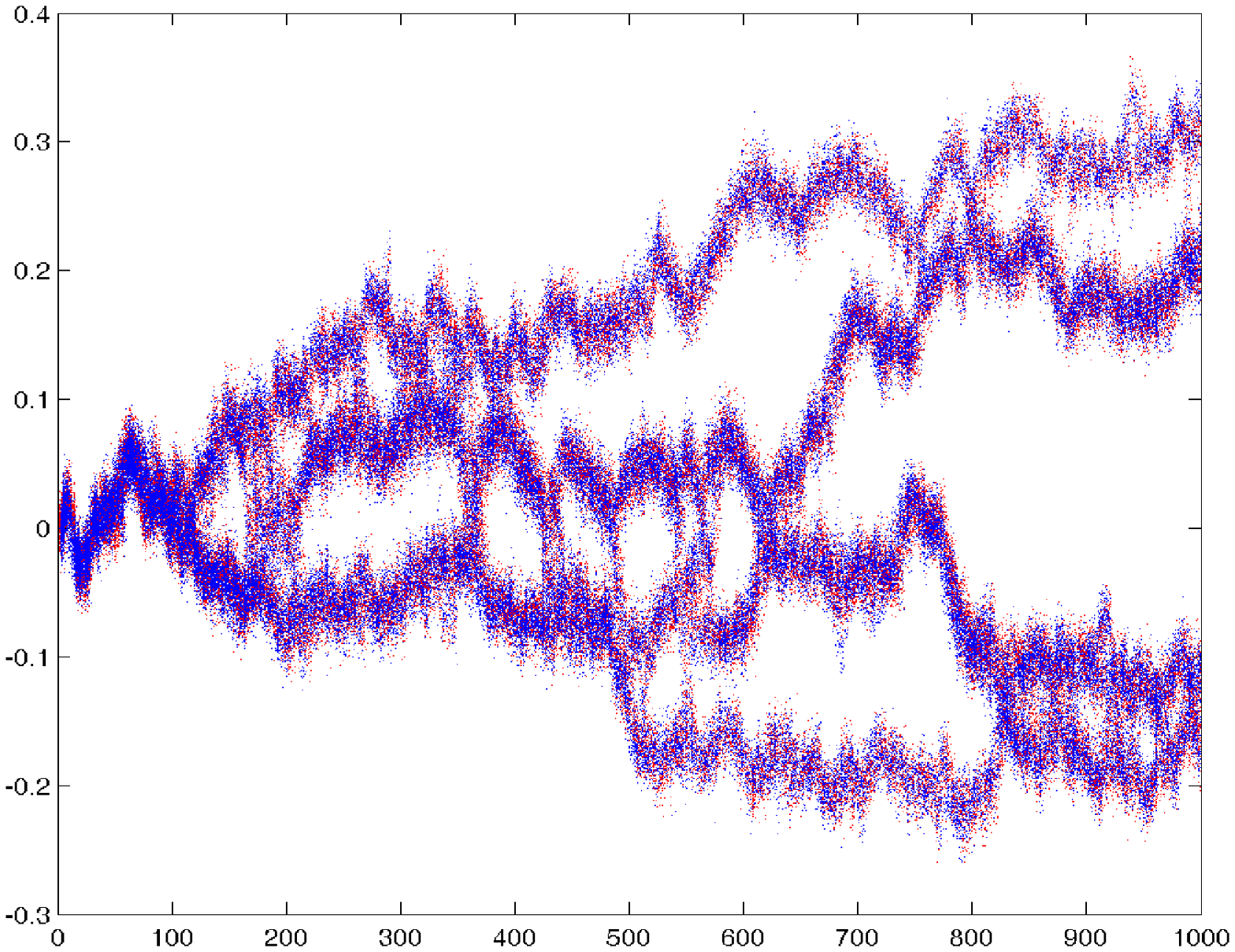}}
  
  \subfloat[]{\label{fig:sim_x3c}%
    \includegraphics[width=0.38\textwidth]{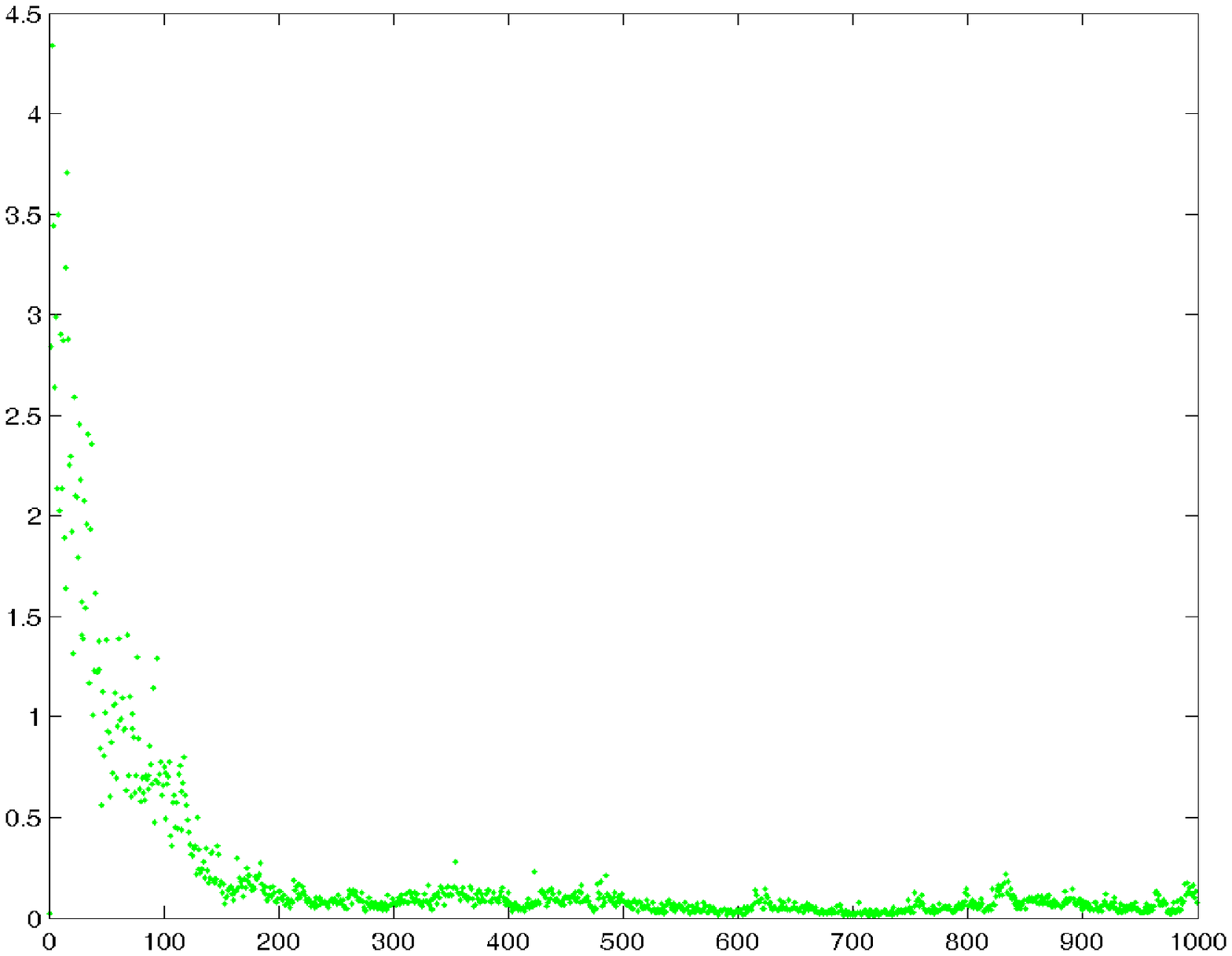}}\hspace{0.1\textwidth}
  \subfloat[]{\label{fig:sim_x3d}%
    \includegraphics[width=0.38\textwidth]{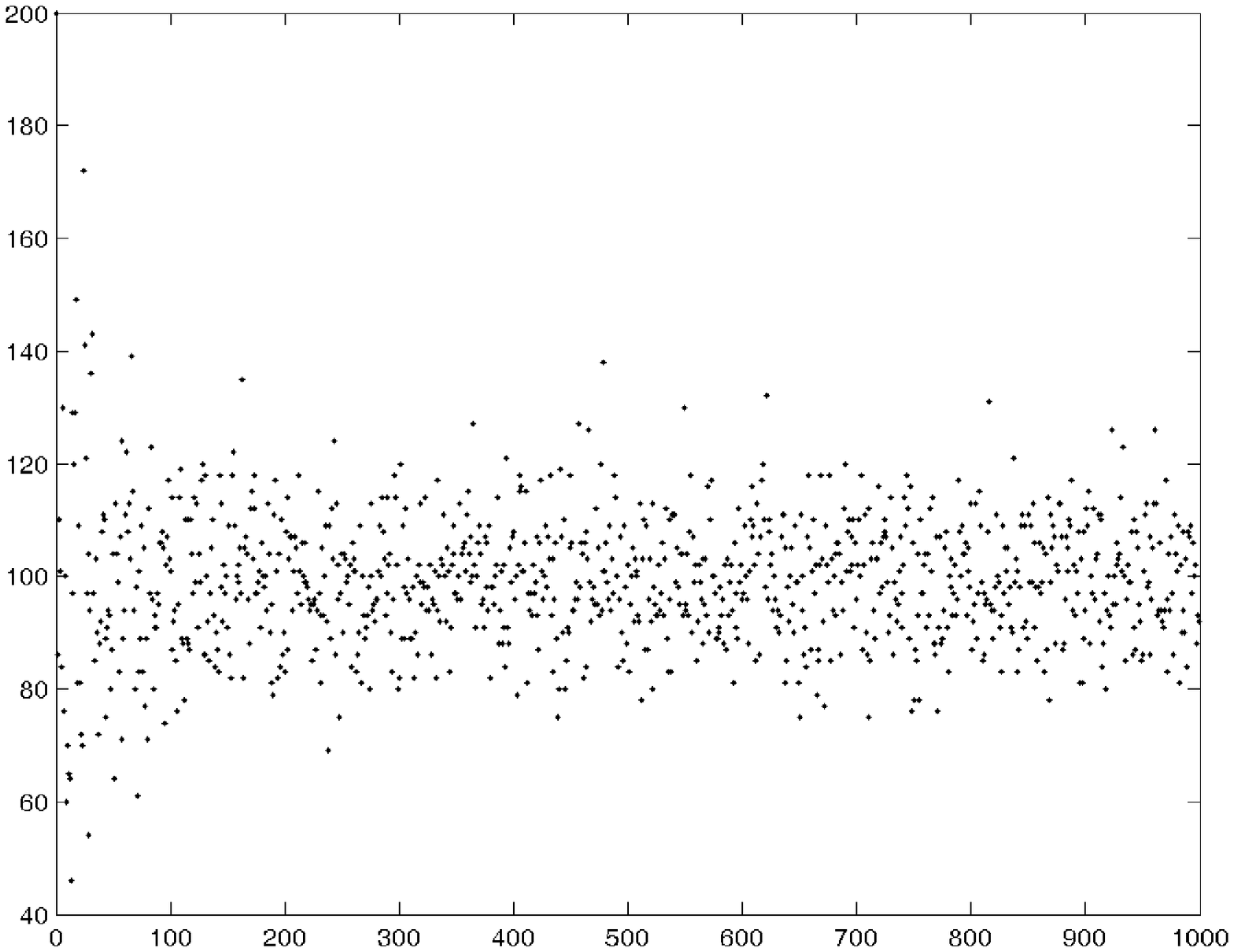}}
  
    \caption{100 individuals, Gaussian food distribution }
  \label{fig:sim_x3}
  \end{figure}

  Some simulation results are shown in Figure~\ref{fig:runAB}, \ref{fig:runC}
  and~\ref{fig:sim_x3}. For the simulations in Figure~\ref{fig:runAB},
  the food resource is concentrated at two points, $x=\pm1$, and the
  population is initially monomorph with phenotype $x=0$. Without
  reinforcement, as in (a), the population remains concentrated around
  $x=0$ although the small mutations are seen as noise in the
  distribution. With reinforcement, as in (b), (c), and (d), the
  population immediately splits in two sub-populations, each one
  exploiting one of the food resources. The graph in (d) shows the
  distribution of the appearance trait $y$. It does also separate in two
  parts, but there is no reason for the parts to stay at any particular
  position, and therefore these will carry out a random walk in the
  $\Y$-space. Eventually the two branches could meet, which would lead
  to the appearance of hybrid phenotypes. The graph in (d) shows the
  evolution of the ``food distribution entropy'',
  \begin{equation*}
    W_c(t) = \sum_{k=1}^{N_t} c_k \log( N_t c_k)\,,
  \end{equation*}
  which is zero if and only if $c_k=1/N_t$ for all $k$. The simulations
  show that the population approaches a situation where all individuals
  attract the same quantity of the food resource.
  
  Figure~\ref{fig:runC} shows the results of a simulation where the food
  resources are equally distributed at the points $x=-1,0,1$, and as
  expected, with reinforcement, the population will then split in three
  sub-populations, but it may happen in different ways: The figures (a),
  (b), and (c), (d) show the results of two simulations with exactly the
  same initial conditions. 
  
  And then  Figure~\ref{fig:sim_x3}, shows the result in which the food
  distribution is Gaussian in $x$, with mean zero, and the $x$-phenotype
  is initially concentrated at $x=2$. We can see in (a) how the whole
  population evolves to $x$-values close to zero, before the splitting
  into sub-populations take place. The graph in (c) shows the food
  distribution entropy, and the graph in (d) the size of the population, $N_t$. 
  
  The model can be reformulated into a more mathematically tractable
  form by identifying the population $Z_t$ by a point measure in $\Z$:
  \begin{equation*}
     Z_t \leftrightarrow \sum_{j=1}^{\langle Z_t,1\rangle}
      \delta_{z_j} \in \M_P(\Z)\qquad  ( N_t = \langle Z_t,1\rangle)\,.
  \end{equation*}
  To find an expression for $Z_{t+1}$, we write the offspring from an
  individual $z_j$ as
  \begin{equation*}
     \Gamma\left(\cdot,z_j\right) = 
            \sum_{i=1}^{\langle
              \Gamma\left(\cdot,z_j\right),1\rangle}\delta_{z_i}
      \qquad ( \mbox{offspring from }  z_j)\,.
  \end{equation*}
  Then the next generation is
  \begin{equation*}
    Z_{t+1} = \sum_{j=1}^{\langle Z_t,1\rangle}
      \Gamma\left(\cdot,z_j\right) = \int_{\Z} \Gamma(\cdot,z) Z_t(dz)\,.
  \end{equation*}
  Here $\Gamma(\cdot,z)$ is a random point measure in $\Z$, whose law
  can be computed from the description above. The details are given
  in~\cite{HenrikssonLundhWennberg2010}. From this we wish to write a
  master equation for the process, to find a formula for 
  \begin{equation*}
    \EE\left[ \int_{\Z} \phi(z) Z_{t+1}(dz)
         \,\Bigg  |\, Z_t 
       \right] = 
       \int_{\Z}{ \EE\left[ \int_{\Z} \phi(z) \Gamma(dz,z')\,\Bigg |\, Z_t
       \right]}Z_t(dz')\,.
  \end{equation*}
  To continue, we first write $\Gamma$ by conditioning on the mating
  partner, 
  \begin{equation}
  \label{eq:matingcond}
   \EE\left[ \int_{\Z} \phi(z)
        \Gamma(dz,z')\,\Bigg |\, Z_t 
      \right] = \int_{\Z} { \EE\left[ \int_{\Z} \phi(z)
          \Gamma(dz,z')\,\Bigg |\, Z_t,z'' 
        \right]} { P(z',z'')Z_t(dz'')}\, , 
  \end{equation}
  where the choice of mating partner, $ \mbox{Prob}(j_k=j)$, is encoded in 
  \begin{equation*}
    P(z',z) = \frac{e^{-|y-y^{*}{}'|^2/ 2\sigma^2}\one_{z\ne z'}}
      {\int_{\Z}e^{-|y''-y^{*}{}'|^2/ 2\sigma^2}\one_{z''\ne z'}Z_t(dz'')}\,.
  \end{equation*}
  The size of the offspring then depends on the resource distribution in
  the population,
  \begin{equation*}
     c(z) = \int_{\X} \frac{e^{-(x-x')^2/2\sigma_x^2}}
              {\int_{\Z}e^{-(x''-x')^2/2\sigma_x^2}Z_t(dz'')}
              \frac{f(dx')}{\|f\|}\,.
  \end{equation*}
  The expectation value in the integral in the
  right-hand side of equation~(\ref{eq:matingcond}) can then be
  computed as as
  \begin{eqnarray*}
    \lefteqn{\EE\left[ \int_{\Z} \phi(z) \Gamma(dz,z')\,\Bigg |\, Z_t,z''
      \right]  =}\\
   &=&  \sum_{k=0}^{\infty}{ \mbox{Prob}[\langle \Gamma(dz,z'),1\rangle = k
      \,|\, Z_t,z''] } \times \\
   & & \qquad \sum_{i=1}^{k} \int_{\Z}
      \phi\left(\frac{z'+z''}{2}+\zeta_i\right)M(\zeta_i)\,d\zeta_i .
  \end{eqnarray*}
  The last sum gives the contribution from each child of parents with
  phenotype $z'$ and $z''$; the phenotype of the child is the average
  of the parents' phenotypes plus a random mutation $\zeta$ which is
  distributed with law $M(\zeta)d\zeta$. The number of offspring $k$ is
  Poisson distributed: 
  \begin{equation*}
    \mbox{Prob}[\langle \Gamma(dz,z'),1\rangle = k
      \,|\, Z_t,z''] = \frac{{\kappa(z',z'')}^k}{k!} e^{-{ \kappa(z',z'')}}\,,
  \end{equation*}
  with 
  \begin{equation*}
     \kappa(z',z'') = \frac{c(z')+c(z'')}{2} \int_{\X}f(dx)\,.
  \end{equation*}
  
  In this form, the equations are still not very explicit, but one may
  formally take the limit of infinitely many individuals as in the paper
  by Méléard and Tran, as described above: we let $\Z_t^n =
  \frac{1}{n}\sum_{i=1}^{n\langle
    Z_t^n,1\rangle}\delta_{(x_i(t),a_i(t))}$, and we {\em assume} that
  this has a limit as $n\rightarrow\infty$, and even that the limit is
  given by a density: $Z_t^n\rightarrow u_t(z)$ with $u_t = u_t(z) =u_t(x,y,y^*) \in L^1(\X\times \Y\times \Y*)= L^1(\Z)$, \quad $u\ge 0$, $
      t\in\mathbb{N}$.
  
  It is then possible to identify the limiting expressions for the food
  distribution, the probability of choosing a particular mating partner,
  et.c.:
  
  \begin{eqnarray*}
        c_k & \mapsto&  c(z;u) = \int_{\X}  \frac{\displaystyle
          e^{-(x-x')^2/2\sigma_x^2}   }
        {\displaystyle \int_{\Z}  e^{-(x''-x')^2/2\sigma_x^2}
          u(z'')\,dz''}
        df(x') \,,\\
        \mbox{Prob}(j_k=j) &\mapsto& \pi(z,z';u)  =
          \frac{ e^{-|y^*-y'|^2  /2\sigma^2}}{ \int_{\Z}
            e^{-|y^*-y''|^2/2\sigma^2} u(z'')\,dz''}\,,\\
         \kappa_k &\mapsto& \kappa(z,z';u) =
        \frac{c(z;u)+c(z';u)}{2}\int_{\X}df(x') \,.     
  \end{eqnarray*}
  
  Finally we may write the master equation for the limiting densities
  $u_n$, 
  \begin{eqnarray}
    \label{eq:3.19}
    \lefteqn{\int_{\Z} \phi(z) u_{t+1}(z)\,dz = } \nonumber \\
    &=&\int_{\Z}\int_{\Z}\int_{\Z}u_t(z')u_t(z'') \lambda(z',z'';u) P_{z'}(z'';u)M(z)
    \phi\left( \frac{z'+z''}{2}+z \right)\, dz'dz''dz \nonumber\\ 
    &=& \int_{\Z} \bigg[
            \int_{\Z}\int_{\Z} u_t(z') u_t(z'')  \lambda(z',z'';u)
       \times \nonumber \\
    & & \hspace{10em}   P_{z'}(z'';u) M\left(
      z-\frac{z'+z''}{2}\right)\,dwdz'' 
  \bigg]\phi(z)\,dz\,.
  \end{eqnarray}
  In the change of variables used to obtain the last line, it is assumed
  that $\Z=\R^d$. To conclude,  $u_{t+1}$ can be expressed in terms of $u_t$
  using the expression in brackets in the last member of equation~(\ref{eq:3.19}).
  However, I want to stress that this is only a very formal derivation,
  and also that there are no mathematical results concerning e.g. the
  long time behavior of the model.
  
  \subsection{An averaging process}
  \label{sec:averaging_process}
  
  If we neglect the rather complicated process of choosing the mating
  partner, the way in which the food resource is distributed, and how
  this affects the number of offspring of a given couple, the process is
  a simple one: according to some probability distribution, choose a
  random couple of individuals, and replace this couple by an offspring
  whose phenotype is the average of the parents' but randomly displaced
  due to mutations. An extremely simplified version of this is the
  following:
  
  Consider $N$ individuals with a scalar phenotype $x\in\R$. The
  population is therefore described by $(x_1,...,x_N)$. The phenotype
  distribution is then updated as follows:
  \begin{itemize}
  \item choose a couple $(x_j,x_k)$ uniformly at random
  \item replace those two individuals with a new couple, as
    \begin{equation*}
       (x_j,x_k) \mapsto \left(\frac{x_j+x_k}{2} +X_1,
        \frac{x_j+x_k}{2} +X_2     \right)\,,
    \end{equation*}
  where $X_1$ and $X_2$ are i.i.d. with probability density $g(x)\in L^1(\R)$.
  \end{itemize}
  
  We now assume that in some limit\footnote{For a given $N$, the
    distribution should be updated until a stationary state has been achieved,
    and then let $N\rightarrow\infty$} all the $x_i$ are distributed with a
  law $f(x)dx$. If $x_0$ is drawn from this distribution, it should be
  the result of a replacement, i.e. 
  \begin{equation*}
    x_0 = \frac{x_{1,1}+x_{1,2}}{2}+ X_0,
  \end{equation*}
  where also $x_{1,1}$ and $x_{1,2}$ are distributed with law $f(x)$,
  and where $X_0$ is a random variable with distribution $g(x)$. The
  same argument can be repeated for  $x_{1,1}$ and $x_{1,2}$, which
  then gives (the notation should be clear)
  \begin{eqnarray*}
     x_0 &=&  \frac{1}{2}\left(
        \frac{x_{2,1}+x_{2,2}}{2}+ X_{1,1} + \frac{x_{2,3}+x_{2,4}}{2}+
        X_{1,2} \right) +X_0\\
     &=&  \frac{1}{4}\left(x_{2,1}+x_{2,2}+x_{2,3}+x_{2,4}+
            x_{2,1}+x_{2,2}\right)  
            + \frac{1}{2}\left( X_{1,1}+ X_{1,2}\right) + X_0\,.
  \end{eqnarray*}
  The procedure can be repeated, and after $n$ iterations we have
  \begin{eqnarray*}
     x_0 &=& \frac{1}{2^n}\sum_{j=1}^{2^n}  x_{n,j} \\
    & & + \;  \frac{1}{2^{n-1}} \sum_{j=1}^{2^{n-1}}
      X_{n-1,j} \;  + \;   ... \; + \\
    & & + \; \frac{1}{4} \sum_{j=1}^{4} X_{2,j} \; + \;
      \frac{1}{2}\left( X_{1,1}+ X_{1,2}\right) \; + \; X_0\,.
  \end{eqnarray*}
  
  By the law of large numbers, the first term converges to $\int_{\R} x
  f(x)\,dx$ when $n\rightarrow\infty$, and the other terms can be
  expressed exactly in terms of the distribution $g$: The law of 
  $ \frac{1}{2^{n-1}} \sum_{j=1}^{2^{n-1}}
  X_{n-1,j} $ is $ 2^{n-1} g^{*2^{n-1}}( 2^{n-1} \cdot)$, for
  example. This gives a relation between the densities $f$ and $g$,
  which is most easily expressed in terms of their Fourier
  transforms:
  \begin{eqnarray*}
    \hat{f}(\xi) &=&  \left(
      \frac{1}{2^n}\hat{f}\left(\frac{\xi}{2^n}\right)^{2^n}
    \right)\times \\
    & &  \left(
      \frac{1}{2^{n-1}}\hat{g}\left(\frac{\xi}{2^{n-1}}\right)^{2^{n-1}}
    \right)\cdots  \left(
      \frac{1}{4}\hat{g}\left(\frac{\xi}{4}\right)^{4}
    \right)  \left(
      \frac{1}{2}\hat{g}\left(\frac{\xi}{2}\right)^{2}
    \right) \hat{g}(\xi)\,.
  \end{eqnarray*}
  If $f$ has bounded second moments, the first factor converges to 1
  when $n\rightarrow\infty$ (this is obviously sufficient to guarantee
  that the law of large numbers holds), and we have then an explicit
  expression of $f$ in terms of the the distribution $g$. One example
  where this can be computed explicitly is when $g$ is a Gaussian
  function with variance $\sigma$; then $f$ will also be Gaussian, but
  with variance $2\sigma$.
  
  We end this section by writing a master equation for the process, and
  computing expressions for the evolution of marginals. While this
  doesn't have much to do with the model for speciation, it gives an
  introduction to much of what will follow.
  For a system of $N$ particles, the configuration space is $\R^N$
  (there is no hardcore condition or similar restriction to the particle
  configuration). Let  $f_N(x_1,.....,x_N)$ be a density in $\R^N$. One
  replacement according to description above, replacing a randomly
  chosen pair of particles by new particles whose position is the
  average of the parents's position plus independent displacements,
  transforms the density as
  \begin{eqnarray}
  \label{eq:master}
  \nonumber    f_N(x_1,.....,x_N) &\mapsto& \; f_N'(x_1,.....,x_N)   = \\
      && \quad \frac{2}{N(N-1)}
  \nonumber  \sum_{j<k} \int_{\R}\int_{\R} f_N(x_1,...,
    x_j',...,x_k',...,x_N)\times \\
     & & \qquad \qquad  g(x_j-\frac{x_j'+x_k'}{2})
    g(x_k-\frac{x_j'+x_k'}{2})\,dx_j'dx_k'
  \end{eqnarray}
  As is commonly done, we assume that the densities are symmetric with
  respect to permutation of the variables, and it is easy to see that if
  this holds for $f_N$, then it also holds for $f_N'$: symmetry is
  preserved by the dynamics. The $k$-particle marginals are defined as 
  \begin{equation*}
     f_{N,k}(x_1,...,x_k) = \int_{\R^{N-k}}
  f_N(x_1,...,x_k,x_{k+1},...,x_N) \, dx_{k+1}\cdots dx_N\,.
  \end{equation*}
  Because of the permutation symmetry, the same result would
  be obtained by leaving any set of $k$ variables untouched, integrating over
  the remaining $N-k$ variables. Integrating both sides of
  equation~(\ref{eq:master}) over $k_{k+1},...,k_N$, we find an
  expression for how the $k$-particle marginals are transformed by the
  replacement process. For $k=1$ and $k=2$ the result is 
  \begin{eqnarray*}
   \lefteqn{ f_{N,1}'(x_1)
  = \left(1-\frac{1}{N}\right)\left(1-\frac{1}{N-1}\right)
  f_{N,1}(x_1)
  +} \\
  & & \qquad\qquad \frac{2}{N}
    \int_{\R}\int_{\R} f_{N,2}(x_1',x_2')
   g(x_1-\frac{x_1'+x_2'}{2})\,dx_1' dx_{2}'\,,
  \end{eqnarray*}
  and
  
  \begin{eqnarray*}
  \lefteqn{  f_{N,2}'(x_1,x_2)  
    = \left(1-\frac{2}{N}\right)\left(1-\frac{2}{N-1}\right)
    f_{N,k}(x_1,x_2)
    + }\\
    & & \qquad\qquad \frac{2(N-2)}{N(N-1)}
    \int_{\R}\int_{\R} f_{N,3}(x_1,x_2',x_{3}')
   g(x_2-\frac{x_2'+x_3'}{2})\,dx_2' dx_{3}'
  + \\
    & &  \qquad\qquad \frac{2(N-2)}{N(N-1)}
    \int_{\R}\int_{\R} f_{N,3}(x_1',x_2,x_{3}')
   g(x_1-\frac{x_1'+x_3'}{2})\,dx_1' dx_{3}'
   + \\
     & &  \qquad\qquad \frac{2}{N(N-1)} \int_{\R}\int_{\R}
  f_{N,2}(x_1',x_2')  g(x_1-\frac{x_1'+x_2'}{2})
  g(x_2-\frac{x_1'+x_2'}{2})\,dx_1' dx_2' \,.
  \end{eqnarray*}
  respectively. Here, and in all higher order terms, we find that the
  expression for $f_k'$ involves terms with $f_{k+1}$, so if $k<N$, the
  system will never be closed. 
  
  The next step is to let $N\rightarrow\infty$. For this to make sense,
  we write $f_{N,k}^j(v_1,...,v_k)$ the distribution obtained after $j$
  replacements, and write
  \begin{equation}
  \label{eq:master2}
    \frac{f_{N,1}^{j+1}(x_1)-f_{N,1}^{j}}{2/N}
  = \int_{\R}\int_{\R} f_{N,2}^{j}(x_1',x_2')
   g(x_1-\frac{x_1'+x_2'}{2})\,dx_1' dx_{2}' - f_{N,1}^{j}(x_1)\,.
  \end{equation}
  Now, if one thinks of $f_{N,k}^{j}$ as being the values  of a time
  dependent function evaluated at discrete points,
  $f_{N,k}^{j}(v_1,...,v_k) =f_{N,k}(v_1,...,v_k,\Delta_t j)$ with
  $\Delta_t= 2/N$, the left-hand side of equation~(\ref{eq:master2}) is a
  finite difference approximation of $\displaystyle \frac{\partial}{\partial
    t}f_{N,k}(v_1,...,v_k,\Delta_t j)$. Passing to the limit as
  $N\rightarrow\infty$ gives
  \begin{equation*}
    \frac{\partial}{\partial t} f_1(x_1,t)  =
    \int_{\R}\int_{\R}
    f_2(x_1', x_2',t)g(x_1-\frac{x_1'+x_2'}{2})\,dx_1' dx_{2}' - f_1(x_1,t)
  \end{equation*}
  If in addition we assume that propagation of chaos holds, that is,
  $f_2(x_1', x_2',t)=f_1(x_1',t) f_1(x_2',t)$ (this will be discussed at
  length in the following sections), then a closed equation for the
  one-particle marginal is obtained:
  \begin{equation}
    \label{eq:beq}
    \frac{\partial}{\partial t} f_1(x_1,t)  =
    \int_{\R}\int_{\R}
     f_1(x_1',t) f_1(x_2',t)g(x_1-\frac{x_1'+x_2'}{2})\,dx_1' dx_{2}' -
     f(x_1,t)\,.
  \end{equation}
  Similar equations can be obtained for all marginals, but if the
  proportion of chaos is assumed to hold, then all information is already
  present in equation~(\ref{eq:beq}).

  \section{Applications from biology: models of flocking animals}
  
  It is fascinating to watch the huge bird flocks flying over big
  cities, or schools of fish that are forming close to bridges,
  sometimes, and recently there have been many attempts to make
  mathematical models to describe the observed phenomena. What is
  intriguing is that these complex structures can be formed without any
  obvious leader, all individuals in the flock should have the same
  status. How can the presumably rather simple rules controlling the
  behavior of individual birds result in this complex collective
  behavior?
  
  In the first section I will present a couple of well known
  mathematical models related to swarming animals (without any claim to
  give a comprehensive list), and then discuss a model which has been
  analyzed in~\cite{CarlenDegondWennberg2011} in some more detail.
  
  \subsection{The Boids and Cucker-Smale models}
  
  The boids model~\cite{Reynolds1987} is a system of ODEs describing the evolution of $N$  particles:
  \begin{eqnarray*}
    \mathbf{\ddot{r}}_i &=& \sum_{j} \Big[
    f(r_{i,j})\mathbf{\hat{r}}_{ij} + \alpha_1\langle
    \mathbf{v}_j-\mathbf{v}_i | r_{ij}<r_c\rangle \\
    && \quad  + \alpha_2 \langle
    \mathbf{r}_j-\mathbf{r}_i | r_{ij}<r_c\rangle
    \Big] -\gamma \mathbf{v}_i + \beta \mathbf{\zeta}_i\,.
  \end{eqnarray*}
  Here $\mathbf{r}_i$ is the position of ``boid'' $j$. The different terms
  describe the boid's desire to move towards the average position of
  swarm, to approach the average velocity of the boids within a smaller
  neighborhood, and to avoid crowding.
  
  There are other models that are discrete in time, and the Cucker-Smale
  model~\cite{CuckerSmale2007} is a particular example for which much
  progress on the 
  mathematical theory has been made recently~\cite{HaTadmor2008}. Here
  the velocities $v_i(t)$, $i=1,...,N$, $t=1,2,3,...$ evolve according
  to 
  \begin{eqnarray*}
    v_i(n+1) -v_i(n) &=& \frac{\gamma}{N}\sum_{j=1}^N a_{ij}\,,
    \left(v_j(n)-v_i(n)\right) \\
    a_{i,j} &=& \frac{1}{(1+\| x_i-x_j\|^2)^{\beta}}\,.
  \end{eqnarray*}
  The main mechanism here is alignment, the particles strive to align
  with the surrounding particles, and the strength of interaction
  depends on the distance between the particles through the function
  $a_{i,j}$. 
  
  A Boltzmann equation inspired by the model of Cucker and Smale has
  been derived by Carillo {\em et
    al.}~\cite{Carilloetal_Flocking2010}. They consider a density of
  individuals, $f(x,v,t)$, interacting pairwise by exchanging velocities
  according to 
  \begin{eqnarray*}
    v^* &=& (1-\gamma a(x-y))v + \gamma a(x-y) w\\
    w^* &=& \gamma a(x-y) v + (1-\gamma a(x-y)) w\,,
  \end{eqnarray*}
  and this leads to a Boltzmann equation
  \begin{eqnarray*}
    \left(\frac{\partial}{\partial t} + v\cdot\nabla_x\right)
    f(x,v,t) &=& Q(f,f)(x,v,t)\,.
  \end{eqnarray*}
  where the collision operator (the binary interactions) is
  \begin{eqnarray*}
    Q(f,f)(x,v) &=& \int_{\mathbb{R}^3}\int_{\mathbb{R}^3}
    \left( \frac{1}{J}f(x,v^*) f(y,w^*)-f(x,v)f(y,w)\right)\, dw\,dy
  \end{eqnarray*}
  
  \subsection{The Vicsec model and a related Boltzmann equation}
  
  A discrete time model somewhat similar to that of Cucker and Smale was
  derived by Vicsec~\cite{Vicsec_etal1995,Czirok1999}, and has since been used in a
  large number of publications. In this model, all velocities have the
  same magnitude, $v_0$, and the direction is updated in one of
  the following ways~\cite{Chateetal2007}:
   \begin{eqnarray}
  \label{eq:vic1}
      v_i(t+\Delta t) = v_0 \theta\left[ \sum_{j\in S_i} v_j(t) + \eta
          \mathcal{N}_i \xi \right]\,,
    \end{eqnarray}
  or
   \begin{eqnarray}
  \label{eq:vic2}
      v_i(t+\Delta t) = v_0 \left( \mathcal{R}_{\eta}\circ\theta\right) \left[
          \sum_{j\in S_i} v_j(t) \right] \,.
    \end{eqnarray}
  In both these models, $S_i$ is the set of neighbors to the particle
  $i$, defined as all other particles in side a ball of given radius
  around $x_i$, the position of particle $i$. This means that particle
  $i$ only interacts with other particles inside this ball. The function
  $\theta[\cdot]$ normalizes a vector: $\theta(v) = v/|v|$. In
  (\ref{eq:vic1}) the new velocity of particle $i$ is computed by first
  taking the average velocity of the particles inside the radius of
  interaction, adding a random vector scaled by the number $N_i$ of
  particles inside the ball of interaction, and finally normalizing the
  magnitude. In (\ref{eq:vic2}) the new velocity is computed by first
  finding the average velocity of the particles in the ball (as in
  (\ref{eq:vic1}), normalizing and finally carrying out a random
  rotation $\mathcal{R}_{\eta}$. The two models have been analyzed
  carefully with respect to e.g. phase transitions
  in~\cite{Chateetal2007}.
  
  A Boltzmann equation related to the Vicsec model has been derived
  in~\cite{BertinDrozGregoire2006}: 
  
  \begin{figure}
    \centering
    \includegraphics{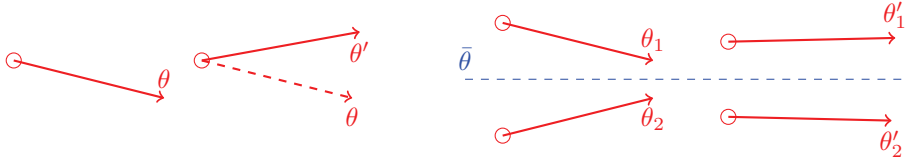}
    \caption{The velocity jumps in the model of Bertin et al.}
    \label{fig:BDGmodel}
  \end{figure}
   \begin{eqnarray*}
  \lefteqn{   \frac{\partial f}{\partial t}(r,\theta,t) + e(\theta)\cdot\nabla
     f(x,\theta,t)  = 
  -\lambda f(r,\theta,t) + \lambda
  \int_{-\pi}^{\pi}d\theta'\int_{-\infty}^{\infty} d\eta p_0(\eta) \times}\\
  &&\hspace*{5cm}\sum_{m=-\infty}^{\infty}
  \delta(\theta'+\eta-\theta+2\pi m) f(r,\theta',t) \\
  &&-f(r,\theta,t)\int_{-\pi}^{\pi} d\theta'|e(\theta')-e(\theta)|
  f(r,\theta,t)\\
  &&+\int_{-\pi}^{\pi} d\theta_1\int_{-\pi}^{\pi}
  d\theta_2\int_{-\infty}^{\infty} d\eta p(\eta)
  |e(\theta_2)-e(\theta_1)|
  f(r,\theta_1,t) f(r,\theta_2,t)\\
  &&\hspace*{5cm}\times \sum_{m=-\infty}^{\infty}
  \delta(\bar{\theta}+\eta-\theta+2\pi m)\,. 
   \end{eqnarray*}

  \subsection{A simple kinetic equation on the circle}
\label{sec:BDG}
  
  In order to derive fluid equations by the Hilbert or Chapman - Enskog
  methods, one needs to know the equilibrium distributions, and in order
  to approach this we will now study a simpler, spatially homogeneous
  model in the plane~\cite{CarlenDegondWennberg2011}, which can be
  derived from a master equation very similar to the one in Section~\ref{sec:averaging_process},
  \begin{eqnarray}
  \label{eq:bdgA}
    \partial_t f(t,\theta) 
  &=&\int_{-\pi}^{\pi}\int_{-\pi}^{\pi}
  \bigg(f(t,\theta')f(t,\theta'+\theta_*) 
  g(\theta-\theta'-\frac{\theta_*}{2}) 
  \nonumber  \\
  &&\quad\qquad - f(t,\theta)f(t,\theta+\theta_*)\bigg)
  \beta(|\sin(\theta_*/2)|)\frac{d\theta'}{2\pi} \frac{d\theta_*}{2\pi}\,,
  \end{eqnarray}
  which corresponds to a jump process as depicted in the rightmost part
  of figure~\ref{fig:BDGmodel}: two particles  with
  velocities $\theta_1$ and $\theta_2$ (so the velocities are
  represented by angles) get new velocities
  $\theta_j'=\bar{\theta_{1,2}}+\theta_j''$, $j=1,2$. In this expression,
  $\theta_1''$ and $\theta_2''$ are two independent angles distributed
  with law $g(\theta)$. Note the similarity with the averaging process
  described in the previous section. The model is also very similar to
  the model of rod alignment in~\cite{Ben-NaimKrapivsky2006}. It is easy
  to see that, for all distributions $g$, the uniform distribution
  $f(\theta)=1/(2\pi)$ is a stationary solution, but it it may not be
  the only one. Equation~(\ref{eq:bdgA}) may be written in terms of the
  Fourier series of $f(t,\theta)$. With
  $f(t,\theta)=\sum_{k=-\infty}^{\infty} a_k(t)e^{ik\theta}$,
   \begin{eqnarray*}
     \frac{da_k}{dt}      &=&\sum_{n} a_{k-n} a_n
     \left(
       \gamma_k \Gamma(n-k/2) -\Gamma(k)
     \right)\,,
   \end{eqnarray*}
  where
   \begin{eqnarray*}
     \gamma_k = \int_{-\pi}^{\pi} g(\theta)\frac{d\theta}{2\pi}
     \qquad\quad
     \Gamma(z) = \frac{\sin(\pi z)}{\pi z}\,.
   \end{eqnarray*}
  ($\Gamma$ actually depends on the function $\beta$ in~(\ref{eq:bdgA}),
  as written here it corresponds to $\beta\equiv1$). To study the linear
  stability of the uniform distribution, we write $f(t,\theta) = 1 + \varepsilon
  \sum_{k=-\infty}^{\infty}b_k(t)e^{ik\theta}$, and then 
  \begin{equation*}
  \displaystyle \frac{d}{dt}b_k(t) =  b_k(t) \underbrace{\bigg( 2 \gamma_k \Gamma(k/2)      -\Gamma(0)-\Gamma(k) \bigg)}_{\lambda_k} + \bigoh(\varepsilon)\,.  
  \end{equation*}
  We then see that up to order $\varepsilon$, the Fourier modes are
  decoupled, and hence the linear stability can established by checking
  the sign of $\lambda_k$. By a direct calculation it follows  that
  $\lambda_k<0$ if 
  $k\ge 2$, and therefore the stability depends on $\lambda_1$, which is
  negative if and only if $\gamma_1<\pi/4$.
  This in turn depends on $g(\theta)$. As an example we take any
  density $\rho(x)\in 
  L^1(\R)$, and let. 
  \begin{equation*}
     g_{\tau}(y) = 2\pi\sum_{j=-\infty}^{\infty}
        \frac{1}{\tau}\rho(\frac{y-2\pi j}{\tau})  \quad\Rightarrow
        \quad \gamma_k = \hat{\rho}(\tau k)\,.
  \end{equation*}
  The parameter $\tau$ determines how concentrated $g$ is around
  $\theta=0$. Clearly, when $\tau\rightarrow 0$, $\gamma_1\rightarrow
  1>\pi/4$, and therefore the first Fourier mode is unstable for
  sufficiently small $\tau$. A similar result can be found
  in~\cite{Ben-NaimKrapivsky2006}. 
  
  \section{Propagation of chaos}
  
  Boltzmann's and Maxwell's kinetic theory was derived from a physical
  point of view, and it would take a very long time before a
  mathematically satisfactory derivation was carried out by
  Lanford~\cite{lanford} for a hard ball gas. And up to date, a derivation
  valid over macroscopic intervals of time is essentially missing (see
  the notes by Pulvirenti for details about
  this~\cite{PulvirentiPortoErcole2010}.
  
  Mark Kac~\cite{Kac1956} invented a Markov process that mimics an
  $N$-particle system, 
  proposed a mathematically rigorous definition of propagation of chaos,
  and showed that his model satisfies this property. In the following
  sections, we will present Kac's model, and his proof, and then follow
  through the steps of e.g. Grünbaum~\cite{Grunbaum} towards an abstract
  theorem stating not only that a large class of Markov processes do
  propagate chaos according to the definition of Kac, but also give
  precise error bounds in terms of the number of particles, and a
  detailed information about the limiting equation. The results are
  proven in~\cite{MiMoWe2011arXiv} and~\cite{MiMo2010arXiv}. An
  important ingredient in the abstract 
  formulation is the de Finetti (or Hewitt Savage) theorem, which is
  also presented in these notes, following the lectures by
  P.L. Lions~\cite{PLL_cours_champsmoyens}.

  \subsection{Kac's approach to the propagation of chaos}

  Kac's model of an $N$-particle system is a jump process on the
   sphere in $\R^N$,
   \begin{equation*}
  \nklot) = \left\{ (v_1,...,v_N) \; \Bigg| \; \frac{1}{2}\left(
      v_1^2 + ... + v_N^2\right) = N\right\}   \,.
   \end{equation*}
  Each coordinate represents the (one dimensional) velocity of a
  particle, and the radius is chosen so that the expected energy of a
  particle (with unit mass) is one. The particles suffer binary
  collisions, which are modeled as jumps involving two coordinates at a
  time: At exponentially distributed time intervals, two coordinates, say
  $v_i$ and $v_j$ are chosen uniformly, randomly, and they are given new
  velocities $v_i'$ and $v_j'$:
  \begin{equation*}
    (v_1,...,v_i,...,v_j,...v_N) \mapsto
     (v_1,...,v_i',...,v_j',...v_N)\,.
  \end{equation*}
  The new velocities are obtained as a random rotation in $\R^2$:
  $\theta$ is chosen at random according to a law $\mu$, and then
  \begin{equation*}
     (v_i',v_j') = ( v_i \cos(\theta) -
          v_j\sin(\theta), v_i \sin(\theta)+v_j \cos(\theta))\,.
  \end{equation*}
  We will use the notation $ V \mapsto V' = R_{i,j}(\theta) V$. Clearly 
  \begin{equation*}
     {v_i'}^2 + {v_j'}^2 = v_i^2 + v_j^2\,,
  \end{equation*}
  and therefore this process preserves energy exactly. On the other hand, there
  is no conservation of momentum
  \begin{equation*}
    {v_i'} + {v_j'} \ne  v_i + v_j \,.
  \end{equation*}
  With only one dimensional velocities only trivial collisions can
  satisfy both energy and momentum conservation.

  The Markov process just described can equivalently be defined by a
  {\em master equation}, which describes the evolution of phase space
  density. Writing $V=(v_1,...,v_N)$, and $\psi(V,t) =
  \psi(v_1,...,v_N,t)$, we assume that the random, initial point of the
  Markov jump process is distributed with density $\psi_0(V)$. To have a
  concrete example, we assume that the law for the random rotations in
   a collision is $\mu(d\theta) = (2\pi)^{-1}d\theta$ (any bounded
   measure can be treated in the same way,  but for singular measures,
   as in e.g.~\cite{SundenWennberg2009}, one needs to be a little more
   careful). The density at time $t$, $\psi(V,t)$ then satisfies
   \begin{equation}
     \label{eq:mastereq}
     \partial_t \psi_N(V,t) = 
     \frac{2}{N-1} \sum_{1 \le
      i < j \le N } \frac{1}{2\pi} \int_{-\pi}^{\pi} \left( \psi_N \big( 
    R_{ij}(\theta)V,t  \big)-\psi_N(V,t)\right)\,d\theta\,.
   \end{equation}
  The factor in front of the sum in the right hand side implies that the
  system jumps on average  $N$ times per unit time, and because the
  coordinates are drawn independently, this means that each coordinate
  is changed approximately twice per unit time. This corresponds to the
  Boltzmann Grad scaling of a real particle system, because each
  particle should then on average  suffer the same (finite) number of
  collisions per unit time, independently of the number $N$ of particles.
  
  Because all particles in a gas are assumed to be identical, the
  probability distribution of initial values should not depend on in
  which order we write them, and this is expressed by saying that the
  initial distribution should be symmetric with respect to permutations:
  
  \begin{definition}
    The density  $ \psi_0(v_1,...,v_N)$ is {\em symmetric} if for any
    pair of variables, $v_i, v_j$,
    \begin{equation*}
      P[  (v_1,...,{ v_i},...,{ v_j},...,v_N) \in A ]= P[
    (v_1,...,{ v_j},...,{ v_i} ,...,v_N) \in A  ]\,.
    \end{equation*}
    The space of configurations obtained by identifying all points 
   $V=(v_1,...,v_N)$ that can be obtained from each other by a
   permutation of the indices,
   \begin{equation*}
       (v_1,..., { v_i},...,{ v_j},..., v_N) \sim
            (v_1,..., { v_j},...,{ v_i},..., v_N)\,,
    \end{equation*}
   is denoted  $\nklot) / \SN$.
  \end{definition}
  Note that the Kac jump process preserves permutation symmetry.

  In the same way as the  Kac master equation corresponds to the
  Liouville equation for a real $N$ particle system, there is a
  Boltzmann equation for the velocity distribution of one particle, that
  can formally be obtained in the limit of infinitely many
  particles. This is {\em Kac's caricature of the Boltzmann equation },
  the Kac equation:
  \begin{equation}
  \label{eq:kacequation}
      \partial_t f(v,t) = \int_{-\infty}^{\infty} \int_{-\pi}^{\pi}
       \left( f(v',t) f(t,w')-f(v,t) f(w,t)\right) \frac{d\theta}{2\pi} dw\,,
  \end{equation}
  where just as  in the definition of the jump process,
  \begin{equation*}
     (v',w') = ( v \cos(\theta) -
          w\sin(\theta), v \sin(\theta)+w \cos(\theta))\,.
  \end{equation*}
  
   In spite of its relative simplicity, its structure is
    very similar to that of the Boltzmann equation, and the two
    equations share many characteristics.
   There are numerous
    studies that consider the trend to equilibrium, tne regularity of solutions,
    its behavior in the presence of external force terms, ... , with
    the hope that it will give insight into the behavior of the the
    full equation. Some relevant references are~\cite{Desvillettes_Kac1995,Fournier2002,BaglandWennbergWondmagegne2007,CarlenGabettaRegazzini2008}.
  
  Mass and energy conservation are among the most important properties
  of the solutions of the Kac equation:
  \begin{eqnarray*}
     \int_{\R} f(v,t)\,dv &=& const\,,\\
     \int_{\R} f(v,t) v^2 \,dv &=& const,
  \end{eqnarray*}
  and also entropy $\int f \log f dv$ is non-increasing, just as for the
  real Boltzmann equation. On the other hand, the momentum is not a
  conserved quantity.

  The master equation and the Kac equation are connected through the
  marginal distributions. We define, for $k=1\cdots N-1$,
  \begin{equation*}
    f^N_k(v_1,v_2,...,v_k,t) = \int_{\Omega_k}
      \psi_N(v_1,...,v_k,v_{k+1},...,v_N,t)d\sigma_k \,,
  \end{equation*}
  where $ \Omega_k = S^{N-1-k}\left(\sqrt{N-v_{1}^2-...-v_k^2}\right)$,
     and $ \sigma_k$ is the uniform normalized measure on $\Omega_k$.
  The marginals $f^N_k$, the ``k-particle distributions'' give the
  distribution of one of the first $k$ coordinates, and because of the
  permutation symmetry, the distribution of any collection of $k$
  different coordinates is the same, and correspond to the joint distribution of
  $k$  randomly chosen particles in a real gas. Because $\psi_N$ is
  assumed to be symmetric, the marginal distributions are too.
  
  The evolution equations for the $k$-particle marginal can be obtained simply by
  integrating the master equation over $v_{k+1}\ldots v_N$. For example,
  \begin{equation}
  \label{eq:Kac1stmarginal}
     \partial_t f^N_1(v_1,t) = \int_{-\sqrt{N-v_1^2}}^{\sqrt{N-v_1^2}}
    \int_{-\pi}^{\pi}\left( f_2^N(v_1',v_2',t) - f_2^N(v_1,v_2,t) \right)
    \frac{d\theta}{2\pi} \,dv_2\,, 
  \end{equation}
  and similar equations can be obtained for all $f^N_k$. If we assume
  that for each $k$, $f^N_k\rightarrow f_k$ for some function
  $f_k(v_1,...,v_k,t)$, which for each $t$ is a density in $\R^k$, then
  it is possible, at least formally,  to pass to the limit
  in~(\ref{eq:Kac1stmarginal}) to get
  \begin{equation}
  \label{eq:Kac1infty}
     \partial_t f_1(v_1,t) = \int_{-\infty}^{\infty}
    \int_{-\pi}^{\pi}\left( f_2(v_1',v_2',t) - f_2(v_1,v_2,t) \right)
    \frac{d\theta}{2\pi} \,dv_2\,.
  \end{equation}
  The situation is similar for all $f_k$: because the right-hand side
  involves $f_{k+1}$, one does not obtain a closed system of
  equations. The whole discussion about propagation of chaos aims at
  proving that if certain hypotheses are satisfied, then
  $f_k(v_1,...,v_k,t) = f_1(v_1,t)\cdot\ldots\cdot f_1(v_k,t)$. The
  interpretation of this is that drawing a $k$-tuple of velocities is
  the same as drawing $k$ velocities independently, 
the particle velocities are
  independent. This cannot be true for any finite $N$, but can sometimes
  be proven to be correct in the limit as $N\rightarrow\infty$.

  \subsection{Propagation of chaos in Kac's model}
\label{sec:Kac}
  
  Kac {\em defined} propagation of chaos as follows:
  
  \begin{definition}
    A sequence of probability measures $
    \psi_N(v_1,...,v_N)$, $ N=1,...,\infty$ is said to have { the
      Boltzmann property}, or to be $ chaotic$ if for each $k$,
    \begin{equation*}
       \lim_{n\rightarrow\infty} f_k^N(v_1,..., v_k) \rightarrow
      \prod_{j=1}^k   \lim_{n\rightarrow\infty} f_1^N(v_j,t)\,.
    \end{equation*}
   Assume that the evolution of a sequence  of probability measures  
    $\psi_N(v_1,...,v_N,t)$ is governed by a family of Markov processes,
    and that the sequence is chaotic for each $t\ge 0$. Then the {\em
      propagation of chaos} is said to hold for these Markov processes.
  \end{definition}
  
  One of the main achievements in~\cite{Kac1956} was Kac's proof that
  propagation of chaos holds for his model, and hence that the Kac
  equation can be derived rigorously as the limit of a many particle
  system.
  
  \begin{thm}[M. Kac]
  Propagation of chaos holds for the master equation~(\ref{eq:mastereq}).
  \end{thm}
  
  Very briefly, the main steps of the proof are as follows:
  Seen as an operator in $L^2(\nklot)$.
    \begin{equation*}
       Q \psi_N (V) = \frac{2}{N-1} \sum_{1 \le
     i < j \le N } \frac{1}{2\pi} \int_{-\pi}^{\pi} \left( \psi_N \big( 
   R_{ij}(\theta)V\big)-\psi_N(V)\right)\,d\theta
    \end{equation*}
  is self adjoint and bounded,and hence $\psi_N(V,t) = \exp( t Q )
  \psi_N(V,0)$, where
  \begin{equation}
  \label{eq:expsum}
   \psi_N(V,t) = \sum_{k=0}^{\infty} \frac{t^k}{k!} Q^k \psi_N(V,0) \,.
  \end{equation}
  
  We then need to compute powers of $Q$. Consider first a bounded
    function $ g_1(V) =   g(v_1)$, {\em i.e.} a function depending on
    only the first component of $V$, and let
    \begin{equation}
  \label{eq:Qg1}
       g_2(V)=g_2(v_1,v_2) = \int_{-\pi}^{\pi} (g_1(v_1\cos(\theta)
       + v_2 \sin(\theta))-g_1(v_1))\frac{d\theta}{2\pi}\,.
    \end{equation}
  By recursion, let
  \begin{eqnarray}
  \label{eq:Qgk}
    \lefteqn{ g_{k+1}(v_1,...,v_k,v_{k+1})  =} \nonumber \\
    &&  \sum_{j=1}^{k}\int_{-\pi}^{\pi}
      ( g_k(v_1,...,v_j\cos(\theta) + v_{k+1}\sin(\theta),...,v_k) -
      g_k(v_1,...,v_k))\frac{d\theta}{2\pi}\,.
  \end{eqnarray}
  The reason for introducing the $g_k$ in this way is that computing 
  $\int_{\nklot} \psi_N(V)g_1(V)$ for all bounded $g_1(V)$ is enough to
  identify the one particle marginal $f^N_1(v_1)$, and that because $Q$
  is self adjoint, $\langle Q\Psi,g_1\rangle =\langle \Psi,Q g_1\rangle
  $. The formulae (\ref{eq:Qg1},\ref{eq:Qgk}) then appear in the calculation of
  $\langle \Psi,Q^k g_1\rangle $.  
  
  Next we assume that the initial data are chaotic, so that
  $\psi_N(v_1,...,v_N,0) = f_{0,k}^N(v_1,...,v_N)$, and that there are
  functions $f_{0,k}$ such that 
  $ f_{0,k}(v_1,...,v_k) = $ \\
  $\lim_{N\rightarrow\infty}
  f_{0,k}^N(v_1,...,v_k)$, and moreover
  \begin{eqnarray*}
  \lefteqn{   \int_{-\infty}^{\infty} ... \int_{-\infty}^{\infty} 
      f_{k+1}(v_1,...,v_{k+1}) g_{k+1}(v_1,...,v_{k+1}) dv_1...dv_{k+1}}\\ 
    &=&\int_{-\infty}^{\infty} ... \int_{-\infty}^{\infty} 
       f_{1}(v_1)\cdots f_1(v_{k+1}) g_{k+1}(v_1,...,v_{k+1}) dv_1...dv_{k+1}\,.
  \end{eqnarray*}
  
  Multiplying all terms in~(\ref{eq:expsum}) with $g_1(v_1)$,
  integrating and letting $N\rightarrow\infty$, we get
  \begin{eqnarray}
  \label{eq:kacseries1}
  \lefteqn{  \int_{-\infty}^{\infty} f_1(v_1,t) g(v_1)\,dv_1 = } \nonumber\\
  &&  \sum_{k=0}^{\infty} \frac{t^k}{k!} 
   \int_{-\infty}^{\infty} ... \int_{-\infty}^{\infty} f_0(v_1)\cdots
   f_0(v_{k+1})
   g_{k+1}(v_1,...,v_{k+1}) dv_1,...,dv_{k+1} 
  \end{eqnarray}
  for $0\le t<2$. Similarly for the two-particle marginals
  \begin{eqnarray}
  \label{eq:kacseries2}
    \lefteqn{ \int_{-\infty}^{\infty}\int_{-\infty}^{\infty}
    f_2(v_1,v_2, t) g(v_1)h_(v_2) \,dv_1dv_2 =} \nonumber \\
  & & \sum_{k=0}^{\infty} \frac{t^k}{k!} 
   \int_{-\infty}^{\infty} ... \int_{-\infty}^{\infty} f_0(v_1)\cdots
   f_0(v_{k+2})
   \gamma_{k+2}(v_1,...,v_{k+2}) dv_1,...,dv_{k+2}\,.
  \end{eqnarray}
  Here the $ \gamma_k$ are obtained by iteration:
  \begin{eqnarray*}
      \gamma_2(v_2,v_2) &=& g(v_1)h(v_2)\\
      \gamma_{k+1} &=&   \sum_{j=1}^k \int_{-\pi}^{\pi}
     (\gamma_k(v_1,...,v_j \cos(\theta) + v_{k+1}\sin(\theta),...,v_k) -
     \gamma_k(v_1,...,v_j,...,v_k))d\theta\,.
  \end{eqnarray*}
  We then need to prove that 
  \begin{equation*}
    \int_{-\infty}^{\infty}\int_{-\infty}^{\infty}
    f_2(v_1,v_2t) g(v_1)h_(v_2) \,dv_1dv_2 =
    \int_{-\infty}^{\infty} f_0(v_1,t) g(v_1)\,dv_1 
      \int_{-\infty}^{\infty} f_0(v_2,t) h(v_2)\,dv_2\,.
  \end{equation*}
  This is done by proving that the series~(\ref{eq:kacseries1})
  and~(\ref{eq:kacseries2}) are convergent, and by comparing the
  terms. This involves expressions like
  \begin{eqnarray*}
    \gamma_3(v_1,v_2,v_3) &=& g_2(v_1,v_3) h_1(v_2) + g_1(v_1)
    h_2(v_2,v_3)\\
   \gamma_4(v_1,v_2,v_3,v_4)& =& g_2(v_1,v_3) h_2(v_2,v_4) +
     g_3(v_1,v_3,v_4)h_1(v_2) \\
    & &+ 
     g_2(v_1,v_4) h_2(v_2,v_3) + g_1(v_1)
     h_3(v_2,v_3, v_4)\,.
  \end{eqnarray*}
  As for the convergence of the series, this turns out to hold of a
  bounded time interval, but this interval can be uniformly estimated,
  and hence one can prove that propagation of chaos holds for any
  bounded time interval.
  
  \subsection{Existence of chaotic states}
  
  Kac proved that there is a large class of functions distributions
  defined as
  \begin{eqnarray*}
    \psi_N(v_1,...,v_N) &=& \frac{  \prod_{j=1}^N   c(v_j)}{ {\displaystyle\int\limits_{\nklot}} \; \prod_{j=1}^N   c(w_j)\;
   d\sigma_N(w_1,...,w_N)} \,.
  \end{eqnarray*}
  The easiest examples are the uniform distributions, which are also the
  equilibria for the Kac master equations. If $\psi_N(V,t)$ is the
  solution to equation~(\ref{eq:mastereq}), then 
  \begin{equation*}
    \lim_{t\rightarrow\infty} \psi_N(V,t)  =  \frac{1}{|\nklot|}\,.
  \end{equation*}
  In this case one can carry out explicit calculations rather easily:
  To compute the limit of a one-particle marginal, let
  \begin{eqnarray*}
    P( v_N> w) =
   \left\{ 
       \begin{array}{ll}
         0 &\mbox{ if } w>\sqrt{N}\\
         \\
         {\displaystyle\int_{v_N>w}} d\sigma_N(v_1,...,v_N) &\mbox{ if }
         -\sqrt{N}<w<\sqrt{N} \\
         \\
         1 & \mbox{ if } w<-\sqrt{N}
       \end{array} \right. 
  \end{eqnarray*}
  Then write the spherical cap as
  \begin{eqnarray*}
     Y_{\sqrt{N},w}&=&  \left\{ (v_{1},v_{2},...,v_{N});    v_{1}^{2}+
     v_{2}^{2}+ ... +   v_{N}^{2} = \sqrt{N}^{2}, v_{N}\ge w   \right\}\,,
  \end{eqnarray*}
  whose  area is 
  \begin{eqnarray*}
      \mu(  Y_{\sqrt{N},w} ) &=& \mu(S^{N-1}(1)) N^{(N-1)/2}
        \int_{0}^{\arccos(w/\sqrt{N})} 
    (1-\cos^{2}(\theta))^{(N-3)/2}\sin{\theta}  \,d\theta \\
   &=& \mu(S^{N-1}(1)) N^{(N-1)/2}  
      \int_{w/\sqrt{N}}^{1}(1-x^{2})^{(N-3)/2}\,dx  \,.
  \end{eqnarray*}
  Finally, because
  \begin{eqnarray*}
      \sqrt{(N-3)/2} \int_{w/\sqrt{N}}^{1}(1-x^{2})^{(N-3)/2}\,dx
    &\rightarrow&
     \int_{w/\sqrt{2}} e^{-x^2}\,dx
  \end{eqnarray*}
  we may deduce that the one-particle distribution converges to a
  Maxwellian, and a similar calculation can be carried out for the
  two-particle distribution, and so on.

  \section{Empirical distributions}
\label{sec:empirical}

  One difficulty with the approach of Kac is that each $N$-particle
  system has its own state space, while the questions of convergence
  would be more easily stated if one could embed all the $N$-particle
  systems in the same space. One approach to this was suggested by
  Grünbaum~\cite{Grunbaum}, who proposed method for proving  a
  propagation of chaos result for 
  the spatially homogeneous Boltzmann equation for hard spheres. We
  begin by discussing this in abstract terms.
  
  The phase space, or configuration space, for the $N$-particles is then
  $E^N$, an $N$-fold product of a Euclidean space $E$, or more precisly,
  in order that the way in which the particles are numbered is not
  important, $ E^N / \SN$, the quotient group of $E^N$ with the
  symmetric group of $N$ elements. That means that 
  \begin{equation*}
    X=(x_1,...,x_N)\in E^N \mbox{  and  }
    \tilde{X}=(\tilde{x}_1,...,\tilde{x}_N)\in E^N 
  \end{equation*}
  are identified if $\tilde{X}$ can be obtained by a permutation of the
  coordinates of $X$.
  We then define the {\em empirical measure} associated with $X$ as
  \begin{equation*}
     \empX = \frac{1}{N}\sum\limits_{j=1}^{N} \delta_{x_j}  \in
         \PP_N(E)\subset \PP(E)\,.
  \end{equation*}
  Here we have introduced the notation $ \PP(E)$ for  the set of
  probability measures on $ E$, and 
       $ \PP_N(E)$ for the  probability measures consisting of $ N$
       Dirac measures of equal mass. This is slightly at variance with
       the usual definition of empirical measure in probability theory,
       where the $x_j$ are assumed to be i.i.d. random variables with
       some distribution $\mu\in\PP(E)$.

  One important property of $ \PP(E)$ is that it is metrizable. Metrics
  can be introduced in several different ways, two of the most commonly used
  metrics being   Lévy-Prokhorov metric, and Wasserstein distance. The
  first of these is defined as follows: Let $ (E,d)$ be a metric space,
  and let $ \PP(E)$ the collection of probability measures on $E$.  
   For $ \mu, \nu\in \PP(E)$,
   \begin{equation*}
      \lpdist(\mu,\nu) =  \inf\left\{
        \varepsilon>0\,|\,\mu(A)<\nu(A^{\varepsilon})+\varepsilon \mbox{
          and } \nu(A)<\mu(A^{\varepsilon})+\varepsilon \mbox{ for all } A\in\BB(E)
   \right\}\,,
   \end{equation*}
  where $A^{\varepsilon} = \{ x\in E\,\mbox{ such that } \inf_{y\in A}
  d(x,y) \le \varepsilon\}$. 
  
  This definition depends, of course, on the
  metric $d$ on $E$. Given a distance on $E$, there is a natural way of 
  introducing a distance on $E^N/\SN$. Let $X=(x_1,...,x_N)$,
  $Y=(y_1,...,y_N)\in E^N/\SN$. Then 
  \begin{equation*}
    \lpdist(X,Y) = \inf_{\sigma\in\SN, \varepsilon>0} 
   \left\{
   \frac{\sharp \left\{ i\; | \; |x_i-y_{\sigma_i}|>\varepsilon
     \right\}}{N}     < \varepsilon
   \right\}\,.
  \end{equation*}
  And with this metric, we may then define the Levy-Prokhorov distance on
  $\MM(E^N/\SN)$. Note that this metric scales well with $N$. For
  example, if $X=(x,x,x,....,x)$ (i.e. $N$ copies of the same $x\in E$)
  and $Y=(0,....,0)$, then $\lpdist(X,Y) = |x|$, independently of
  $N$. On the other hand, if $X=(x_1,0,...,0)$, then we always have
  $\lpdist(X,Y)\le 1/N$.  
  
  As for the Wasserstein distance, it is defined as follows: 
  Let   $ \Gamma(\mu,\nu)$ be the collection of $ \gamma\in\PP(E\times E)$ 
   such that 
   $ \mu = \int_{E} \gamma(\cdot,dy)$ and
   $ \nu = \int_{E} \gamma(dx,\cdot)$. Hence the $\gamma$s are are joint
   probability measures with $\mu$ and $\mu$ as marginal distributions,
   and the Wasserstein distance is defined as 
   \begin{equation*}
      W_p(\mu,\nu)^p = \inf\limits_{\gamma\in\Gamma(\mu,\nu)} \int_{E\times
     E} d(x,y)^p\,d\gamma(x,y)
     =\inf \Ee\left[ |X-Y|^p\right]\,.
   \end{equation*}
  
  \subsection{The Hewitt-Savage theorem}
  
  Hewitt-Savage theorem (which is an extension of  de Finetti's
  theorem) that is the topic 
  of this section, is relevant for the discussion of propagation of
  chaos, but it is included here also  to serve as an
  introduction to the rather abstract notation that will be used later. 
  The material  is essentially taken from~\cite{PLL_cours_champsmoyens}.
  
   Let $\mu^N\in\PP(E^N/\SN)$.
   The {\em  marginals}
     $ \mu_k^N\in \PP(E^k/\SN)$ are defined by 
     \begin{equation*}
    \int_{E^k} \phi(x_1,...,x_k) \mu_k^N(dx_1,...,dx_k) =    
      \int_{E^N} \phi(x_1,...,x_k) \mu^N(dx_1,...,dx_k,...dx_N)\,,
     \end{equation*}
  which should hold for every  symmetric $ \phi\in \CC(E^k)$. In the
  following we assume that $E$ is a compact metric space, but for
  example $E=\R^n$ could be treated in a similar manner, although with
  some extra technical complication.
  
  We then identify $X\in E^N/\SN$ with an empirical measure as above,
  \begin{equation*}
     \hat{\mu}_X=\frac{1}{N}\sum_{k=1}^N \delta_{x_k} \in \PP_N(E)\subset \PP(E)\,,
  \end{equation*}
  this then yields a natural identification of  functions  $ u: E^N/\SN
  \rightarrow \R$ with functions  $ \  \PP_N(E)\rightarrow \R$:
  \begin{equation*}
     u(X) \leftrightarrow u\left(\frac{1}{N}\sum_{j=1}^N
           \delta_{x_j}\right)\,.
  \end{equation*}
  Note that if this identification does not automatically preserve
  properties like continuity, unless som care is taken in chosing the
  metric on $\PP(E^k/\SN)$. The Levy-Prokhorov example given above has
  this property.
  
  Next we consider a sequence of functions $\left\{ u_N:  E^N/\SN\rightarrow
       \R\right\}_{N=1}^{\infty}$. These are all defined on different
     spaces, and hence there is no immediate way of comparing the
     functions, and talking about convergence, et.c.. But the
     identification with measures in $\PP(E)$ provides a mean of doing
     so: one can say that the sequence $u_N$ converges if the
     corresponding sequence of measures converges.
  
  We have the follwing compactness result:
  
     Consider a bounded sequence  $ \left\{ u_N :  E^N/\SN \rightarrow
       \R \right\}$ ,
   with $ |u_N|\le C$, and let $\omega:
   \R^+\rightarrow\R^+$ be a strictly decreasing function with
   $\omega(r)\xrightarrow[r\rightarrow0]{} 0$.  Assume that 
   \begin{equation*}
      \left| u_N(X)-u_N(Y) \right|  \le \omega( \lpdist(X,Y) )\,.
   \end{equation*}
  What this says is that the sequence is bounded and uniformly
  continuous with modulus of continuity $\omega$.
   Then there is a subsequence $ u_{N'}$  and $ U\in
   \CC(\PP(E))$ such that 
   \begin{equation*}
     \left\| u_{N'}(x_1,...,x_{N'}) - U\left(\frac{1}{N'}\sum_{j=1}^{N'}
     \delta_{x_j}\right)
   \right\|\rightarrow 0\,.
   \end{equation*}
   Of course, in most
     cases the limiting measure cannot be identified with any function
     $u_N:  E^N/\SN\rightarrow \R$ for a finite value of $N$, but would
     are in all cases  good examples of  symmetric functions of
     infinitely many  variables $x_i\in E$. 
  
  Next consider the following calculation: with $ E$ compact, let
     $ \mu^N\in\PP(E^N/\SN)$, $ N=1,2,3....$, and consider the
     marginal distributions $ \mu^N_k = \int
     \mu^N(\cdot,...,dx_{k+1},...,dx_N)\in \PP(E^k / \SSS_k)$, for $1\le
     k < N$. Because
     $E$ is compact, $\PP(E^k/\SSS_k)$ is also compact, and for every 
   $k$ there is a subsequence $ (N')$ such that 
   \begin{equation*}
     \mu^{N'}_k \rightarrow \mu_k\in\PP(E^k/\SSS_k)\,.
   \end{equation*}
  By the usual diagonal procedure, it is then possible to extract a
  subsequence $ N''$  such that 
  \begin{equation*}
     \mu^{N''}_k \rightarrow \mu_k\in\PP(E^k/\SSS_k)\;\; \mbox{\color{black}
       for all } \;\; k \,.
  \end{equation*}
  By construction, the $ (\mu_k)_{k=1}^{\infty}$ satisfy
  \begin{equation}
  \label{eq:measurecompatibility}
        \mu_k = \int_{E}  \mu_{k+1}(\cdot,dx_{k+1})\,.
  \end{equation}
  The Hewitt-Savage theorem~\cite{HewittSavage1955}, which is a
  generalization of a theorem by de Finetti, concerns sequences of
  measures with exactly this property. It is common to express this as
   exchangability: a sequence of random variables, $X_1,X_2,....$
  is said to be {\em exchangeable} if for any $n\ge1$ and any
  $\sigma=(\sigma_1,...,\sigma_n) \in\SSS_n$, the $n$-tuples
  $X_1,...,X_n$ and $X_{\sigma_1},...,X_{\sigma_n}$ have the same distribution. 
  
  \begin{thm}
  \label{thm:hewittsavage}
    Assume that a sequence of measures  $ \mu_k\in \PP(E^k/\SN),\;\;
    k=1,2...$ satisfies (\ref{eq:measurecompatibility}). There is $
    \pi\in\PP(\PP(E))$ such that for all $ k\ge 1$, 
   \begin{equation*}
   \mu_k(dx_1,...,dx_k) = \mathlarger{\int} \prod_{i=1}^k  m(dx_i) \,  \pi(dm)\,. 
  \end{equation*} 
  
  \end{thm}

  The easiest conceivable example, $
  \pi = \delta_{\bar{m}},\quad  \bar{m}\in\PP(E)$. Then $\pi$ is  a measure on
  $\PP(\PP(E))$ that is concentrated on $\bar{m}\in\PP(E)$, and
  \begin{equation*}
   \mu_k(dx_1,...,dx_k) = \mathlarger{\int} \prod_{i=1}^k  m(dx_i) \,
   \delta_{\bar{m}}(dm) =  \bar{m}(dx_1)\cdots\bar{m}(dx_k) \,.
  \end{equation*}
  That is: if $ \pi$ is concentrated in one point $ \bar{m}$, the measures
   derived from $ \pi$ factorize.
   
  \medskip
  
  \noindent
  {\sc Proof of theorem~\ref{thm:hewittsavage}
    (P.L. Lions)}~\cite{PLL_cours_champsmoyens}: With $E$ compact,
  $\PP(E)$ is a compact metric space, and we have constructed functions
  $U\in C(\PP(E))$. One can also define polynomials on $\PP(E)$, and
  these obviously also belong to $C(\PP(E))$. The constants, polynomials
  of degree zero, are in  $C(\PP(E))$, and to define 
  polynomials of degree 1, take  $ \phi\in\CC(E)$, and let
  \begin{equation*}
   m \mapsto \int_E \phi(x) m(dx) \equiv P_1(m)\,.  
  \end{equation*}
  If $m_1\ne m_2$, one can find $\phi$ such that $\int_E \phi(x) m_1(dx)
  \ne \int_E \phi(x) m_2(dx)$, and hence these linear functions separate
  points in $\PP(E)$. Then the monomials $ P_j(m)$ of degree $j$ are defined as
  follows. Take $ \phi_j\in\CC(E^j/\SN)$ and then let
  \begin{equation*}
    P_j(m) =  \int_{E^j} \phi_j(x_1,...,x_j) m(dx_1)\cdots m(dx_j)\,.
  \end{equation*}
  From these definitions one may then define polynomials of all orders,
  and  the Stone-Weierstrass theorem states that the set of
  polynomials is dense in $C(\PP(E))$. Evaluating these polynomials on
  empirical measures  $ m= \frac{1}{k} \sum_{j=1}^k \delta_{x_j}$, we
  find 
  \begin{eqnarray*}
      P_j(m) &=& \int_{E^j} \phi_j(y_1,...,y_j) \left(\frac{1}{k}
        \sum_{i=1}^k \delta_{dx_i}(dy_1)\right)\cdots\left(\frac{1}{k}
        \sum_{i=1}^k \delta_{dx_i}(dy_j)\right) \\
      &=& \frac{1}{k^j}\sum_{i_1=1}^k\cdots \sum_{i_j=1}^k
   \phi_j(x_{i_1},...,x_{i_j}) +..
   \end{eqnarray*}
   For example,
   \begin{equation*}
     P_2\left(\frac{1}{3}(\delta_{x_1} +\delta_{x_2}
       +\delta_{x_3})\right) =  \frac{1}{9}\sum_{i,j=1}^3 \phi(x_i,x_j)\,.
   \end{equation*}
  In these second degree polynomials there are  some terms like
  $\frac{1}{9}\phi(x_1,x_1)$, and the same will happen for polynomials
  of higher degree. However, when $k$ is much larger than the
  degree $j$ of the polynomial, a vast majority of the terms will be of the
  form $\phi(x_{i_1},...,x_{i_j})$ where all the arguments are different. 
  
  Now let  $ \{\mu_k\}_{k=1}^{\infty}$ be the measures given
  in the statement of the theorem, and consider
  \begin{equation*}
       \int P_j\left(\frac{1}{k} \sum_{j=1}^k \delta_{x_j}\right)
  \mu_k(dx_1,...,dx_k) \approx
   \int \phi(x_1,...,x_j) \mu_k(dx_1,...,dx_j,...dx_k) \,.
  \end{equation*}
  The difference between the left and right terms is due to the presence
  of terms with two or more of the arguments of $\phi$ are taken to be
  the same $x_i$, and so vanishes when $k$ is large compared to $j$. The
  error can be estimated by a simple combinatorial argument. And 
  \begin{equation*}
    \int \phi(x_1,...,x_j) \mu_k(dx_1,...,dx_j,...dx_k) = 
   \int \phi(x_1,...,x_j) \mu_j(dx_1,...,dx_j)
  \end{equation*}
  because of the relation~(\ref{eq:measurecompatibility}). 
  
  Next we {\em define} a linear functional on the set of polynomials  $
  P\in\CC(\PP(E))$ by
  setting
  \begin{equation*} 
    \ell(P) =
   \lim_{k\rightarrow\infty}  \int_{E^k} P\left(\sum_{j=1}^k \delta_{x_j}
   \right) \mu_k(dx_1,...,dx_k)
  \end{equation*}
  Then $P\mapsto\ell(P)$ is  positive and  $\ell(1)=1$, so $\ell$ is a
  positive, bounded functional defined on a dense subset of
  $\CC(\PP(E))$, and can be extended to all of   $\CC(\PP(E))$. Then
  Riesz's theorem states that there is a measure  $ \pi\in \PP(\PP(E))$ such
  that
  \begin{equation*}
    \ell(U) = \int_{\PP(E)} U(m) \; \pi(dm)\,. 
  \end{equation*}
  This measure, $\pi$ is the the desired measure. We only need to check
  that the measures  $\mu_k$ can be obtained from $\pi$ as stated in
  the theorem.  
  To this end, consider  
  \begin{equation*}
    \int_{\PP(E)} \prod_{i=1}^k m(dx_i) \,
  \pi(dm) \in \PP(E^k/\SN) \,.
  \end{equation*}
 Integrating a function   $ \phi(x_1,...,x_k)\in \CC(E^k/\SN)$
  with respect to  this measure gives
  \begin{eqnarray*}
      \int_{E^k} \phi(x_1,...,x_k) \int_{\PP(E)} \prod_{i=1}^k m(dx_i) \,
       \pi(dm) &=&\int_{\PP(E)} { P_k(m)}\, \pi(dm)  = \ell(P_k)  \\
  &= &    \int_{E^k} \phi(x_1,...,x_k)
    \mu_k(dx_1,...,dx_k) \,.
  \end{eqnarray*}
  and this completes the proof.
  \qed
  
  \medskip
  
  By a short calculation we can now establish that   if   $ \pi$ gives
  rise to measures that factorize, then $ \pi$ is a Dirac measure:
  
  \begin{prop}
   Assume that $ \pi\in\PP(\PP(E))$, and that 
   \begin{equation*}
      \int_{\PP(E)} m(dx_1) m(dx_2) \pi(dm) = \left(\int_{\PP(E)} m(dx_1) \pi(dm)\right)
   \left( \int_{\PP(E)}  m(dx_2) \pi(dm)\right) \,.  
   \end{equation*}
  Then there is $ \bar{m}\in\PP(E)$ such that $ \pi = \delta_{\bar{m}}$
  \end{prop}
  
  \noindent
  {\sc Proof:}(Lions)~\cite{PLL_cours_champsmoyens}. 
    Multiply by $ \phi(x_1)\phi(x_2)$, $
   \phi\in\CC(E)$, and integrate. Then the lefthand side is
   \begin{equation*}
          \int_{\PP(E)} \left(\int_E \phi(x_1) m(dx_1)\right)
            \left(\int_E \phi(x_2) m(dx_2)\right) \pi(dm) =
     \int_{\PP(E)} \left(\int_E \phi(x) m(dx)\right)^2 \pi(dm)\,.
   \end{equation*}
   and the right hand side
   \begin{equation*}
     \left(\int_{\PP(E)}  \left( \int_E \phi(x)  m(dx)\right)  \pi(dm)\right)^2\,.
   \end{equation*}
   But Jensen's inequality states that 
   \begin{equation*}
      \left(\int_{\PP(E)}  \left( \int_E \phi(x)  m(dx)\right)  \pi(dm)\right)^2
    \le
     \int_{\PP(E)} \left(\int_E \phi(x) m(dx)\right)^2 \pi(dm)\,,
   \end{equation*}
   with equality only if $ \int_E \phi(x) m(dx) = c(\phi)$, i.e.
   independent of $ m$ on the support of $ \pi(dm)$. It follows that
   $ \pi = \delta_{\bar{m}}$.
   \qed

  \section{Estimates on the  propagation of chaos for $N$-particle
      systems}
  
  Kac's  approach to the propagation of chaos concerns a
  very simplified model, and it is not a trivial matter to extend it to
  more realistic models. For example, his model is essentially
  Maxwellian, which means that the collision rate of two particles does
  not depend on their relative velocity. Grünbaum~\cite{Grunbaum}
  circumvented some of these problems by a more abstract approach based
  on identifying an $N$-particle configuration $(v_1,...,v_N)\in
  \R^{3N}$ with an empirical measure $\frac{1}{N}\sum_{j=1}^N
  \delta_{v_j}$, very much like in the discussion about the
  Hewitt-Savage theorem. Some other works in the same direction are the
  results on statistical solutions to the Boltzmann equation that can be
  found e.g. in~\cite{ArkCapIan1991}. 
  
  In Grünbaum's terminology, the set $\PP(E^N/\SSS_N)$ is convex, and
  its extreme points are exactly the symmetric Dirac measures, 
  $\frac{1}{|\SSS_N|}\sum \delta_{(v_{\sigma_1},...,v_{\sigma_N})}$. Any
  point $\beta\in\PP(E^N/\SSS_N)$ can be expressed as the barycenter of
  the extreme points: there is a measure $\Omega_{\beta}$ such that
  $\beta = \int_{E^N/\SSS_N} X \Omega_{\beta}(dX)$, and his work is
  based on an analysis of the evolution of $\Omega_{\beta}$ under the
  collision process.
  
  The remaining part of this paper is a summary of the results
  in~\cite{MiMoWe2011arXiv} and in~\cite{MiMo2010arXiv}, where the
  method of Grünbaum is rephrased, and new quantitative results on the
  rate at which the propagation of chaos is achieved with an increasing
  number of particles.

  \subsection{The abstract setting} 
\begin{figure}[h!]
  \centering
\includegraphics{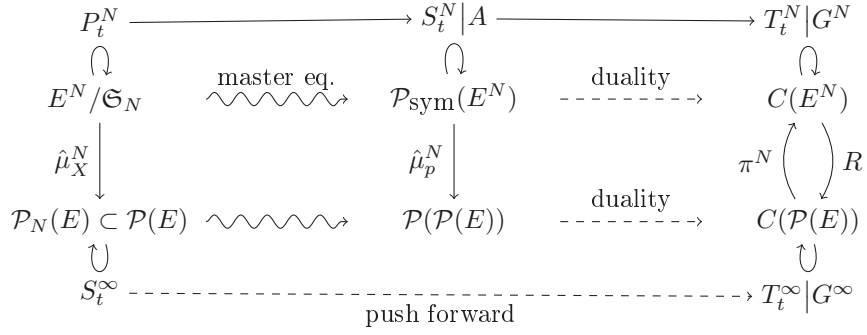}
  
  \caption{A summary of spaces and their relations. Semigroups are in
    most cases given together with their generators, as in $S^N_t \big
  | A$.}
  \label{fig:diagram}
\end{figure}

To formulate the main results
  in~\cite{MiMoWe2011arXiv} and in~\cite{MiMo2010arXiv}, we need to
  introduce a number of spaces, operators on the spaces and maps
  bestween the spaces, as shown in Figure~\ref{fig:diagram}
\medskip

\medskip
   
    We consider a family of Markov jump processes, $
      N=1,2,3,....$, on the spaces $ E^N/\SN$. The space  $E$
      is a locally compact, separable metric space. For every  $N$
      there is a process $X_t = (x_{1,t},...,x_{N,t})$ and a propagator so that
      \begin{equation*}
        X_t = P_t^{N}X_0\,,
      \end{equation*}
      and because we want to be able to identify $X_t$ and $\tilde{X}_t$
      if the components of $\tilde{X}_t$ can be obtained as a
      permutation of the components of $X$, we ask that $P_t^N$ commutes
      with permutations of the components. This is the {\em
        microscopic} description, each componenent $(x_i)_t\in E$
      representing the position of one particle in the $N$-particle system.

  Two most elementary examples  are the Kac model, and
  Grünbaum's model of the three dimensional Boltzmann equation. The
  techniques developed here works also in many other cases, and some
  more examples are given later in this paper. In the diagram in
  Figure~\ref{fig:diagram}, the phase space, and the propagator
  $P_t^N$ are shown in the upper left corner of the diagram.

  The Markov processes can be descibed by the master equation (or
  Kolmogorov equation): 
      Let $ p^N_t=\loi(X_t)$, i.e.,  
      \begin{equation*}
        P[ X_t\in A\subset E] =
      \int_A p^N_t(dx)\,.
      \end{equation*}
  There is a { semigroup} $ S^N_t: \PPs(E^N) \rightarrow
      \PPs(E^N)$ such that $ p^N_t = S^N_t p^N_0$. 
      This semigroup has a { generator}  $ A$, so that $ p_t$
      satisfies
      \begin{equation*}
        \partial_t p^N_t = A p^N_t\,.
      \end{equation*}
This is represented in the middle, upper part of the diagram.
      There is also a dual semigroup  $ T^N_t:  \CC(E^N/\SN)\rightarrow \CC(E^N/\SN)$
     with a corresponding generator  $ G^N$, shown in the upper right
     part. The two semigroups are
     related as follows: 
      For all $ p^N\in\PP(E^N)$, $ \phi\in\CC(E^N)$,
      \begin{equation*}
        \langle p^N, T^N_t(\phi)\rangle = \langle S^N_t(p^N),\phi\rangle\,,
      \end{equation*}
      and, with $ \phi_t = T^M_t \phi$,
      \begin{equation*}
        \partial_t \phi_t = G^N \phi_t, \quad \phi_0=\phi\,.
      \end{equation*}

Thus the upper part of the diagram represent the $N$-particle system
in three different ways, essentially equivalent.
  In kinetic theory we are intersted in rigorously deriving the {\em
    Boltzmann equation} as a limit of an $N$-particle system, and in
  Kac's work~\cite{Kac1956}, this corresponds to derving the nonlinear
  Kac equation from his $N$-particle model. In this abstract setting we
  assume that there is a formal {\em mean field description}, and
  equation that governs the evolution of a one-particle distribution,
  $p_t\in\PP(E)$:
  \begin{equation}
  \label{eq:boltzmann}
     \partial_t p_t = Q(p_t)\,,
  \end{equation}
  and typically this is a nonlinear equation. For the purpose of this
  paper, we require that the initial value problem to
  equation~(\ref{eq:boltzmann}) has a unique solution for initial data
  $p_0 \in \PP(E)$ or some subset of $\PP(E)$. The solution is
  represented by a semigroup, $p_t = S^{\infty}_t p_0$. We see this in the
  lower left of the diagram. The lower part of the diagram thus
  concerns the limit as $N\rightarrow\infty$ in the $N$-particle
  system, and it also provides the arena for comparing the solutions
  to the $N$-particle system and the Boltzmann equation that is the
  formal limit. 
And the objective, is to
  prove that, given that certain conditions are satisfied, the
  one-particle marginals converge to the solution of
  equation~(\ref{eq:boltzmann}):
  \begin{equation*}
     p^N_{1,t} = \int_{E^{N-1}}
    p^N_t(\cdot, dx_2,...,dx_N) \xrightarrow[N\rightarrow\infty]{}  p_t.
  \end{equation*}
  In order to proceed with this, we first need to represent an
  $N$-particle configuration $X_t=(x_{1,t},....,x_{N,t}$ in $\PP(E)$
  in the lower part of the diagram. This representation is provided by
the map $\hat{\mu}_{\cdot}^{N}$, which takes a point $X$ as an
argument, and returns a point measure:
\begin{equation*}
 X\mapsto \hat{\mu}_X^N= \sum_{j=1}^N \delta_{x_j}\in \PP(E)_N\subset
\PP(E) \,.
\end{equation*}

If $X$ is random, distributed according
  to a law  $ \mathcal{L}(X) = p_t^N\in\PP_{\mbox{sym}}(E^N)$, 
the resulting measure $\hat{\mu}_{X}^{N}$ is random with a
distribution  which, in the
diagram, is denoted  $\hat{\mu}_{p}^N \in \PP(\PP(E))$, as indicated in the middle
    column. 

In the same way that $\PP_{\mbox{sym}}(EÊ)$ is related by duality to
the set of continuous, symmetric functions, here denoted $C(E^N)$,
there is a duality relation between $\PP(\PP(E))$ and
$C(\PP(E))$. Clearly the exact properties of this duality depends
strongly on the topology on $\PP(E)$.

The maps $\pi^N$ and $R$ between $C(E^N)$ and $C(\PP(E))$ are
  defined through $\hat{\mu}_X^N$ as follows: For $\Psi\in C(\PP(E))$,
  \begin{equation*}
  \pi^N: \Psi \mapsto \psi(x_1,...,x_N) =
  \Psi(\hat{\mu}_{x_1,...,x_N}^N)\,.    
  \end{equation*}
 That is, given $X=(x_1,...,x_N)$,
  the argument of $\psi$, we get a measure
  $\hat{\mu}_{x_1,...,x_N}^N\in \PP_N(E)$, and this measure is then
  taken as an argument when evaluating  $\Psi$. 

Conversely, given $\phi\in C(E^N)$, the function $\Psi = R \phi \in
C(\PP(E))$ is defined as
\begin{equation*}
\PP(E) \ni \mu \mapsto \int_{E^N}
\phi(x_1,.,,,.x_N) \mu(dx_1)\mu(dx_2)\cdot\ldots\cdot\mu(dx_N)\,.  
\end{equation*}
In
the terminology of Section~\ref{sec:empirical}, $R\phi$ is a monomial
of degree $N$.

The last objects in the diagram are $T^{\infty}_t$ and
  $G^{\infty}$. The former is the push forward of $S^{\infty}_t$. For
  $\Psi\in C(\PP(E))$,  $T^{\infty}_t \Psi$ is defined by
  $\PP(E)\ni\mu\mapsto \Psi( S^{\infty}_t \mu)$, and $G^{\infty}$ is
  its generator. Note that here $T^{\infty}_t$ is a {\em linear} semigroup.

The relation between the non-linear semigroup $S^{\infty}_t$ and the
linear semigroup $T^{\infty}_t$ is simular to the relation between a
$Hamiltonian$, finite dimensional dynamical system, and the
corresponding Liouville eqation. Consider a (deterministic) system of
ODEs  in $ \R^n$, ({\em e.g.} Hamiltonian): 
\begin{eqnarray}
\label{eq:dynsyst}
   \left\{ 
     \begin{array}{l}
    \dot{x} = F(x)\\
         x(0) = x_0
     \end{array}\right.\,,
  \end{eqnarray}
with solution $ x(t) = S^t x(0)$. The Liouville equation states how a phase space density
  $ \Phi_t(x)$ is transported by the flow $ F$:
  \begin{align*}
 \Phi_t(x) = \Phi_0(S^{-t}x) := {T'}^t\Psi_0(x)   \,.
  \end{align*}
Here $\Phi_t(x)$ is given explicitly in terms of $\Phi_0$, but the
expression involves $S^{-t}$, and hence is not valid if
equation~(\ref{eq:dynsyst}) cannot be solved backwards. The remedy is
to study the dual problem:  Take $ \varphi\in C(R^n)$, multiply and integrate:
  \begin{equation*}
 \int_{\R^n}    \Phi_0(\underbrace{S^{-t}x}_{ y}) \varphi(x) \,dx =  \int_{\R^n}
 \Phi_0(y) \varphi(S^t y) \,dy =  \int_{\R^n}
 \Phi_0(y) T^t \varphi( y) \,dy\,.  
  \end{equation*}
The linear semigroup $T^t$ is here defined through the forward
evolution of $x$. In Figure~\ref{fig:diagram}, the Boltzmann
equation is the deterministic dynamical system, but in the phase space
$\PP(E)$. A phase space density in $\PP(\PP(E))$ is transported by the
flow via $S^{\infty}_{-t}$, but in general $S^{\infty}_t$ is not
reversible, and therefore it may be only the dual representation that
makes sense.

Solutions to the equation 
\begin{equation*}
  \partial_t \Psi = G^{\infty} \Psi
\end{equation*}
in $C(\PP(E))$ are known as {\em statistical  solutions} to the
Boltzmann equation, and have been studied  for example
  in~\cite{ArkCapIan1991}.

\subsection{The main result, and important hypotheses on spaces and operators}

Like in Section~\ref{sec:Kac}, the $N$-particle system is represented
by a family of master equations, one for each $N$. That means that for
each $N$ we consider
\begin{equation*}
  p_t^N = S^N_t p^N_0\,,
\end{equation*}
where $ p_0^N\in \PP(E^N/\SN)$. The formal limit as
$N\rightarrow\infty$ is given by the Boltzmann equation, whose
solution is 
\begin{equation*}
  p_t = S^{\infty}_t p_0\,,
\end{equation*}
where $p_0\in \PP(E)$.

\begin{thm}(Very informally)
\label{thm:main}
 There is a  constant $ C(k,\ell) >0$ only
 depending on $ k$
   and $ \ell$ such that for any
   $  \varphi = \varphi_1 \otimes \varphi_2 \otimes \dots \otimes \,
   \varphi_\ell $
   with  $ N \ge 2 \ell$:
   \begin{equation}
\label{eq:theoremestimation}
     \sup_{[0,T]}\left| \left \langle
         \left( S^N_t(p_0 ^N) - \left( S^\infty _t (p_0) \right)^{\otimes N}
         \right), \varphi  \right\rangle
     \right|  \le C(k,\ell,N) \rightarrow 0
     \;\mbox{\color{black}when}\; N\rightarrow\infty     
   \end{equation}
\end{thm}

What this says is that if we compare the solution to the $N$-particle
master equation with an $N$-fold tensor product of the solution to the
limiting Boltzmann equation, only through the distribution of the
first $\ell$ particles, the difference decreases as
$N\rightarrow\infty$, and we can compute the rate  explicitly.

Obviously the statement cannot be true in this generality. To begin
with, we must of course make very precise the statement that the
Boltzmann equation is the formal limit of the $N$-particle
system, both in terms of the equations and in terms of the initial
data. The proof can be seen as perturbation result, where the 
$N$-particle systems are treated as perturbations of the limiting
equation, and because of that, the nonlinear semigroup $S^{\infty}_t$
must satisfy a rather strong regularity condition.

In addition to this, because the actual estimates are carried out in
the framework indicated by the Figure~\ref{fig:diagram}, 
and for everything to
work, we must be very precise when defining the spaces. In particular,
the test-functions $\varphi$ in equation~(\ref{eq:theoremestimation})
must be taken from $\in (\FF_1 \cap \FF_2 \cap \FF_3)^{\otimes \ell}$,
where the $\FF_j$ are subspaces of $C(E)$ which are defined below.

So here are then the four main hypotheses on the abstract semigroups and the
spaces they are acting on. While seemingly complicated, they can be
readily verified in some relevant cases, and and example of this will be given
later.

  \begin{itemize}
  \item[{\sc (H1)}] {\sc Convergence of the generators:}{\em  There exists some
  integer $ k \ge 1$,  $ \theta \in (0,1]$ and a space
  $ \FF_1 \subset C(E)$ such that
  \begin{displaymath}
    \forall \, \Phi \in C^{k,\theta}(\FF'_1), \quad 
    \left\| \left( G^N \, \pi^N - \pi^N \, G^\infty \right)  \,  
      \Phi \right\|_{L^\infty(E^N)} \le \epsilon(N) \, 
     \| \Phi \|_{C^{k,\theta}(\FF_1')}
  \end{displaymath} 
  for some function $ \epsilon(N)$ going to $ 0$ as
  $ N$ goes to infinity. }
\noindent
Here $ \FF_1'$ is the dual of $ \FF_1$, and $\PP(E)\subset \FF_1'$,
and $ C^{k,\theta}(\FF_1')\subset C(\PP(E),\FF_1')$ denotes the set of 
  Hölder differentiable functions on $ \PP(E)$. This must be defined,
  of course.  $ \CC(\PP(E),\FF_1')$ is the set of contiuous functions on
  $ \PP(E)$ defined with the topology given by $ \FF_1'$.

 \item[{\sc (H2)}] {\sc Differential stability of the limit
     semigroup:} {\em  We
  assume that for some affine space $ \FF_2 \subset \FF_1$ such that
  $ \FF'_2$ is a Banach space, the flow $ S_t ^\infty$ on
  $ P(E)$ is
  $ C^{k,\theta}(\FF'_1,\FF'_2)$ uniformly on $ [0,T]$,
  for some integer $ k >0$
  and $ \theta \in (0,1)$: there exists $ C_T ^\infty >0$ such that
{
  \begin{displaymath}
    \sup_{[0,T]}   \| S^\infty _t \|_{C^{k,\theta}(\FF'_1,\FF'_2)} \le
    C_T ^\infty.
  \end{displaymath}
}
}
   This implies, for example, that the associated pushforward semigroup
  $ T^\infty _t$ maps $ C^{k,\theta}(\FF_2 ')$ into  $
  C^{k,\theta}(\FF_1')$. Also this hypothesis relies on a stringent
  defintion of Hölder differentialbility in these spaces.

\item[{\sc (H3)}] {\sc Weak stability of the limit semigroup:} 
  {\em There is a space  $ \FF_3 \subset C(E)$,
  such that }
{
  \begin{displaymath}
    \forall  \, \mu, \nu \in P(E)  \qquad 
    \sup_{[0,T]} \mbox{dist}_{\FF'_3} 
    \left( S^\infty _t(\mu), S^\infty _t(\nu) \right) \le
    C_T^{\infty,w} \, \mbox{dist}_{\FF'_3} (\mu,\nu).
  \end{displaymath}
}
In other words,  $ T_t ^\infty$ propagates the
$ C^{0,1}(\FF'_3)$ norm. 

\item[{\sc (H4)}] {\sc Compatibility of the projection:} {\em  We assume that 
  the dual of $ \FF_3$,  $ \FF_3 '$ satisfies}:
{  \begin{displaymath}
    \left\| \pi^N (\Phi) \right\|_{\FF_3 ^N} =
  \left\| \Phi \circ \hat \mu^N _X \right\|_{\FF_3 ^N} 
  \le C_\pi \, \left\| \Phi \right\|_{C^{0,1}(\FF_3 ')}.
\end{displaymath}
}
Hence the space $\FF_3$ and its dual are defined so as to give the
maps between $C(E^N)$ and $C(\PP(E))$ as shown  to the right in
Figure~\ref{fig:diagram} good properties, and this is also why {\sc
  (H3)} is needed in addition to {\sc (H2)}.
\end{itemize}

With the definitions implicitly given in these hypotheses, it is
possible to express the constant in
Equation~(\ref{eq:theoremestimation}) in more detail:

\begin{multline*} C(k,\ell, N) = C(k,\ell) \, \Bigg[
    \frac{\|\varphi\|_{L^\infty(E^\ell)}}{N} + T \, C_T ^\infty \, \epsilon(N)
    \, \|\varphi\|_{\FF_2 ^k \otimes (L^\infty)^{\ell-k}} \\
    + C_\pi \, C_T ^{\infty,w} \, \mbox{{\em dist}}_{\FF'_3} \left(p_0 ^N,
      p_0 ^{\otimes N}\right) \, \| \varphi \|_{\FF_3 \otimes
      (L^\infty)^{\ell-1}} + C_T
    ^{\infty,w} \, \Omega_N ^{\FF'_3} (p_0) \, \, \| \varphi \|_{\FF_3
      \otimes (L^\infty)^{\ell-1}}\Bigg]\,,
  \end{multline*}

  where 
  \begin{displaymath}
      \Omega_N ^{\FF'_3} (p_0) := \int_{E^N} \mbox{{\em dist}}_{\FF'_3} 
    \left(\hat \mu^N _X, p_0 \right) \, 
    p_0 ^{\otimes  N}(dX)\,.
  \end{displaymath}
 
The first term is related to the converegence of the generators (as in
{\sc (H1)}, the second term to {\sc (H2)}. The third term simply says
that the initial data to the $N$-particle system must be close to a
tensor product, and the last term, finally, that the initial data
$p_0$ must be well approximated by an empirical distribution.

This means that if
  $  \mbox{{\em dist}}_{\FF'_3} \left(p_0 ^N, p_0 ^{\otimes N}\right) 
     \xrightarrow[N \to \infty]{} 0$ and 
   $  \Omega_N ^{\FF'_3} (p_0) \xrightarrow[N \to \infty]{} 0$
  propagation of chaos holds, with an explicitly computable rate
  depending on $ \epsilon(N)$, $ \mbox{{\em
      dist}}_{\FF'_3} \left(p_0 ^N, p_0 
    ^{\otimes N}\right)$ and $ \Omega_N ^{\FF'_3} (p_0)$.

\subsection{Differential calculus on $ \PP(E)$.}

The requirements of the semigroups, as given in {\sc(H1)} and
{\sc(H2)} above are expressed in terms of differentiability of functions
functions $ \Phi : \PP(E)\rightarrow \R$
and semigrioups $ S_t^{\infty}: \PP(E)\rightarrow  \PP(E)$. The exact
definitions are given in this section, together with a couple of examples.

\begin{definition}
  Let  $ \GG_1$ be an affine metric space
 and $ \GG_2$ Banach space, and let  $ \MM^j(\GG_1,\GG_2)$ be the set of
  bounded $j$-multilinear maps  from   $ \GG_1$ to $ \GG_2$. We say
  that  $ \psi: \GG_1\rightarrow \GG_2$ belongs to  $
C^{k,\theta}(\GG_1,\GG_2)$, the space of functions
  $ k$ times
  differentiable with $ \theta$ H\"older regularity from $
  \GG_1$ to $ \GG_2$,  if there exist $
D^j\psi:\GG_1\rightarrow\MM^j(\GG_1,\GG_2)$ such that
{
  \begin{displaymath}
    \forall \, \mu, \, \nu \in \GG_1 \qquad 
    \left\| \psi(\nu) -  
    \ \sum_{j=0}^k \left \langle D^j \psi(\mu) ,  
    (\nu - \mu)^{\otimes j} \right \rangle \right\|_{\GG_2}  
  \le C \, \mbox{{\em dist}}_{\GG_1} (\mu, \nu)^{k+\theta}.
\end{displaymath}
}
The norm is
\begin{align*}
\| \psi \|_{C^{k,\theta}(\GG_1,\GG_2)} &= 
 \qquad \sum_{j=1} ^k \left\| D^j \psi
 \right\|_{C(\GG_1,\MM^j(\GG_1,\GG_2))}  \\
&\qquad + 
  \sup_{\mu,\nu \in \GG_1} \frac{ \Big\| \psi(\nu) -
  \sum_{j=0}^k \langle D^j \psi(\mu) , (\nu - \mu)^{\otimes j}
  \rangle \Big\|_{\GG_2} }{\mbox{{\em dist}}_{\GG_1}(\mu,\nu)^{k+\theta}}\,.
\end{align*}
\end{definition}
In this paper, $\GG_2$ is either $\R$ or  a subset of
  $ \CC(\PP(E))$.
Note that for $ \theta =0$, continuity is not required: $ C^{0,0} =L^\infty$.

\bigskip

As a first example we show that polynomials are differentiable.
Take $ \FF = \mbox{Lip}(E)$, and $ \FF' = (\PP(E),d_{Lip})$, where the
Lipschitz distance is given by $ d_{Lip}(\mu,\nu) =
\sup_{\phi\in\mbox{Lip}(E)}\left\{ \int_E\phi(x)
  (\mu(dx)-\nu(dx))\;:\; \|\phi\|_{\mbox{Lip}}=1\right\} $. A
monomial  of degree $ k$ in $ \CC(\PP(E))$ is defined by
  \begin{displaymath}
 \mu \mapsto P_k(\mu) = \int_{E^k} \psi(x_1,...,x_k)\mu(dx_1)\cdots\mu(dx_k)\,,   
  \end{displaymath}
where $\psi\in Lip(\E^k)$.  We first compute $ P_k(\nu)=P_k( \mu + (\nu-\mu))$:
  \begin{eqnarray*}
 P_k(\nu) &=&
 \int_{E^k}\psi(x_1,...,x_k)\underbrace{(\mu_1+(\nu_1-\mu_1)) \cdots (\mu_k+(\nu_k-\mu_k)) }_{
   \mu_j:=\mu(dx_j)}  \\
&= & P_k(\mu)\\
&  &+ { \underbrace{ \sum_{j=1}^k \int_{E^k} \psi(x_1,...,x_k)
        \mu_1\cdots\mu_{j-1}\mu_{j+1}\cdots\mu_k (\nu_j-\mu_j)}_{
        T_1}}\\
&  &+ { \underbrace{ \sum_{i<j} \int_{E^k}\psi(x_1,...,x_k)
     [\mu_1\cdots\mu_k]_{i,j} (\nu_i-\mu_i)(\nu_j-\mu_j)}_{T_2}}\\
   &  & + \cdots
  \end{eqnarray*}
Here $ [\mu_1\cdots\mu_k]_{i,j} =
\mu_1\cdots\mu_{i-1}\mu_{i+1}\cdots\mu_{j-1}\mu_{j+1}\cdots\mu_k$, and 
 $ T_1$ represents the first term in a Taylor
  expansion, $ T_2$ the second etc.  The first term, $ T_1$ can be rewritten
  \begin{align*}
     \sum_{j=1}^k& \int_{E^k} \psi(x_1,...,x_k)
        \mu_1\cdots\mu_{j-1}\mu_{j+1}\cdots\mu_k (\nu_j-\mu_j) =\\
    &\int_{E}  \underbrace{\sum_{j=1}^k \int_{E^{k-1}} \psi(x_1,...,
        \underbrace{y}_{\mbox{pos} j},...,x_{k-1}) \mu_1\cdots\mu_{k-1}
      }_{P_{k-1}(\mu;\cdot)}(\nu(dy)-\mu(dy) )\,,
  \end{align*}
where $P_{k-1}(\mu,y)$ is a polynomial in $\mu$, of degree $k-1$,
parameterized by $y$, 
and $P_{k-1}(\mu;\cdot) \in \mbox{Lip}(E)$. As a function of  $y\in  E$ it
is Lipschitz continuous and  $ P_{k-1} (\mu;\cdot)\in \MM^1(\PP(E),
\R)$ by duality. Finally
  \begin{displaymath}
 |P_k(\nu)-P_k(\mu)|\le \|P_{k-1}(\mu;\cdot)\|_{\mbox{Lip}}\;\; d_{\mbox{Lip}}(\mu,\nu)\,.
  \end{displaymath}
Therefore these polynomials are once differentiable, but as with
polynomials in $\R^n$, the calculations yield polynomials of a lower
degree, and therefore it is possible to differentiate again.

The second example is directly related to the propagation of chaos
estimates, and shows that  $ T^{\infty}_t$  is differentiable in $t$.
Take $ \Phi\in\CC^{1,\theta}(\FF')$ and $
p_0\in\FF'$. Then, by definition
\begin{eqnarray*} 
\left( G^\infty \Phi \right) (p_0)&:=& {d \over dt} (T^\infty_t
  \Phi) (p_0) |_{t=0}\,,
\end{eqnarray*}
and, from the diagram, $(T^\infty_t
  \Phi) (p_0) = \Phi( S^{\infty}_tp_0 ) = \Phi( p_t)$.
Therefore 
\begin{eqnarray*} 
\left( G^\infty \Phi \right) (p_0)&=&  {d \over dt} \Phi (p_t) |_{t=0}
= \lim_{t \to 0}   { \Phi(p_{t}) - \Phi(p_0)
    \over t} \\
  &=& \lim_{t \to 0} \left\{ \left\langle D\Phi[p_0], {p_{t} - p_0 \over
        t} \right\rangle + \OO \left( { \mbox{dist}_{\FF'}( p_{t},
        p_0)^{1+\theta} \over
        t} \right) \right\} \\
    &=& \left\langle D\Phi[p_0], {d p_t \over dt} |_{t=0} \right\rangle =
    {\left\langle D\Phi[p_0], Q(p_0) \right\rangle}.  
\end{eqnarray*}
Here we have used the definition of differentiability for $\Psi$, and
the arrive at a formula for $G^{\infty}$ in terms of $Q$ the generator
of the nonlinear semigroup $S^{\infty}$.

\subsection{Proof of the abstract theorem}
 
The purpose of this section is to prove the
estimate~(\ref{eq:theoremestimation}), but many of the details are
left out, as they can be found in~\cite{MiMoWe2011arXiv} and in~\cite{MiMo2010arXiv}.

 Take $ \varphi\in
 (\FF_1\cap\FF_2\cap\FF_3)^{\otimes\ell}$. Then~(\ref{eq:theoremestimation})
 can be split in several terms as  
 \begin{eqnarray*}
   &&\left| \left \langle \left( S^N_t(p_0^{\otimes N}) - \left( S _t ^\infty
           (p_0)
         \right)^{\otimes N} \right), 
       \varphi 
       \otimes {\bf 1}^{\otimes N-\ell} \right\rangle \right| \le
   \hspace{\fill}  
   \\ \\
   &&\le \left| \left\langle S^N_t(p_0^N), 
      \varphi  \otimes
      {\bf 1}^{\otimes N-\ell} \right\rangle -
     \left \langle S^N_t(p_0^N), 
       R[\varphi] \circ \hat\mu^N_X \right\rangle  \right|
   \\ \\
   &&+ \left| \left\langle p_0^N, T^N_t ( R[\varphi] \circ
       \hat\mu^N_X) \right\rangle
     - \left\langle p_0^N, 
       (T_t ^\infty R[\varphi] ) \circ \hat\mu^N_X) \right\rangle  \right|
   \\  \\
   &&+ \left| \left\langle p_0^N, 
       (T_t ^\infty R[\varphi] ) \circ \hat\mu^N_X) \right\rangle
     - \left\langle p_0^{\otimes N}, 
       (T_t ^\infty R[\varphi] ) \circ \hat\mu^N_X) \right\rangle  \right|
   \\  \\
   &&+ \left| \left\langle p_0^{\otimes N}, (T_t ^\infty R[\varphi] ) \circ
       \hat\mu^N_X) \right\rangle - \left\langle (S _t ^\infty (p_0))^{\otimes \ell} ,
       \varphi \right\rangle \right| =: { \TT_1 + \TT_2 + \TT_ 3 + \TT_4}.
 \end{eqnarray*} 
Of these terms, 
 $ \TT_1$ is controlled by purely combinatorial arguments, but the
 other terms depend on the hypotheses stated above. Thus the
 the consistency estimate {\sc (H1)} on
   the generators plus the fine stability assumption {\sc (H2)} on the limit
   semigroup gives an estimate of  $ \TT_2$,
and  the $ \TT_3$, involving the  chaoticity of the initial data
depends on  measure stability assumption {\sc (H3)} on the limit
   semigroup, and the compatibility condition {\sc (H4)} on
   $\pi^N$. Finally, 
  $ \TT_4$ is controlled in terms of the function $ \Omega_N ^{\FF_3
     '}(p_0)$ (measuring how well $p_0$ can be approximated in weak $ \FF_3
   '$ distance by empirical measures), and also this estimate relies
   on the weak measure   stability assumption {\sc (H3)}.

\medskip
\noindent
{\sc Estimate of $ \TT_1$}: For this term,
 \begin{displaymath}
   \TT_1 := \left| \left\langle S^N_t(p_0^N), 
       \varphi \otimes
       {\bf 1}^{\otimes N-\ell} \right\rangle - \left \langle S^N_t(p_0^N),
       R[\varphi] \circ \hat\mu^N_X \right\rangle \right|\,.
 \end{displaymath} 
Here $ \varphi \otimes {\bf 1}^{\otimes N-\ell}$ is the function
$(x_1,...,x_N)\mapsto \phi(x_1,...,.x_{\ell}$, and because of the
symmetry of  $S^N_t(p_0^N)$ it can be replaced by the symmetrized
version $ \left(\varphi 
       \otimes {\bf 1}^{\otimes N-\ell} \right)_{sym} $, which is
     obtaind as a normalized sum over all permuations of the variables
     $x_1,...,x_N$. Also,
 \begin{displaymath}
 R[\varphi] \circ \hat\mu^N_X  = \pi^N R[\phi] = \int_{E^{\ell}}
 \phi(y_1,...,y_{\ell})   \hat{\mu}_X(y_1)\cdots 
   \hat{\mu}_X(y_{\ell})\,,
 \end{displaymath}
and therfore, by estimating the number of terms
with  $ \ell$ different coordinates $ x_i$ in $ \phi$, we find
 \begin{displaymath}
   \forall \, N \ge 2 \ell, \qquad  
   \left| \left(\varphi 
       \otimes {\bf 1}^{\otimes N-\ell} \right)_{sym} -
     \pi^N R[\varphi] \right| 
   \le 2 \, \frac{\ell^2 \, \| \varphi \|_{L^\infty(E^\ell)}}{N}\,,
 \end{displaymath}
 for any $ \varphi \in C_b(E^\ell)$. It is essentially the same
 calculation as in the proof of the Hewitt-Savage theorem in
 Section~\ref{sec:empirical}.

\medskip
\noindent
{\sc Estimate of $ \TT_2$}: Here we wish to prove that, for
  any $ t \ge 0$ and any $ N \ge 2 \ell$
 \begin{multline*}\quad
   \TT_ 2  :=  \left|  \left\langle p_0^N, 
     T^N_t ( R[\varphi] \circ \hat\mu^N_X) \right\rangle 
     - \left\langle p_0^N, \left(T_t ^\infty R[\varphi] ) \circ
       \hat\mu^N_X \right) \right\rangle \right| \\
   \le  C(k,\ell) \, C_T ^\infty \, \epsilon(N) \, \| \varphi \|_{\FF_1 ^k
     \otimes (L^\infty)^{\ell-k}}\,,
 \end{multline*}
 where   $ C(k,\ell)$ is a constant depending only on  $ k$ and $ \ell$.
 The proof is based on the following calculation, in which the
 role of the generators of the semigroups is made visible:
 \begin{equation*}
 T^N_t \pi^N - \pi^N T^\infty _t = - \int_0 ^t \frac{d}{ds} \left( T^N
   _{t-s} \, \pi^N \, T^\infty _s \right) \, ds = \int_0 ^t T^N _{t-s} \,
 \left[ G^N \pi^N - \pi^N G^\infty \right] \, T^\infty _s \, ds.  
 \end{equation*}
 From the hypothesis {\sc (H1)} it follows  that, for any $ t \in [0,T]$
 \begin{eqnarray*}\nonumber
   &&\left\| (  T^N_t \pi^N R[\varphi] 
     - \pi^N T^\infty _t R[\varphi] )  \right\|_{L^\infty(E^N)}  \\
   &&\qquad \le T  \, \epsilon(N) \, \sup_{s \in [0,T]}  
   \left\| T^\infty_s R[\varphi]\right\|_{C^{k,\theta}(\FF_1 ')}.
 \end{eqnarray*}

The next step is to estimate  $ T^{\infty}_t(R[\varphi])=
R[\varphi](S^{\infty}_t \cdot) \in C^{k,\theta}(\FF'_1)$,
and that computation is carried out using the differential calculus
developed above, and the fact the chain rule applies also here. The
details can be found in~\cite{MiMoWe2011arXiv}.
\medskip

\noindent
{\sc Estimate of $ \TT_3$}
Take $ t \ge 0$ and $ N \ge 2 \ell$. The desired estimate here is
 \begin{align*}
 \TT_3 =& \left| \left\langle  \left( 
  p_0^N - p_0 ^{\otimes N} \right), \left(T_t^\infty R[\varphi]
  \right) \circ \hat\mu^N_X 
         \right\rangle \right| \\
  &\qquad\qquad   \le \mbox{dist}_{\FF' _3} \left(p_0 ^N,p_0 ^{\otimes N}\right) \, 
  \left\|  \left(T_t ^\infty R[\varphi] \right) \circ \hat\mu^N_X
  \right\|_{\FF_3 ^N}\,.
 \end{align*}
 But using first {\sc (H4)} and then {\sc(H3)} gives
 \begin{eqnarray*}
  \left\| \left(T_t ^\infty R[\varphi] \right) \circ \hat\mu^N_X
 \right\|_{\FF_3 ^{\otimes N}} &=& \left\|\pi^N \left(T_t ^\infty R[\varphi] \right)
 \right\|_{\FF_3 ^{\otimes N}}\\
&\le& C_{\pi} \, C_T ^{\infty, w} \, \left\| R[\varphi] \right\|_{C^{0,1}(\FF'_3)}\,.
\end{eqnarray*}
and the calculation is completed with
   $ \varphi = \varphi_1 \otimes \cdots \otimes \varphi_\ell 
     \in \FF^{\otimes
     \ell}$, $k \in \N$, $\theta \in (0,1]$.  
\medskip
\noindent
{\sc Estimate of $ \TT_4$}:
Here we need to estimate 
 \begin{displaymath}
   \TT_4 := 
   \left| \left\langle p_0^{\otimes N}, 
     \left(T_t ^\infty R[\varphi] \right) \circ \hat\mu^N_X \right\rangle
     - \left\langle \left(S _t ^\infty (p_0)\right)^{\otimes \ell} , 
                 \varphi \right\rangle \right| \equiv 
    \bigg|  \TT_{4,1} -  \TT_{4,2} \bigg|\,.
 \end{displaymath}
 for $ t\ge 0$ and $ N\ge 2 \ell$.
  The first term can be written
 \begin{eqnarray*}
   \TT_{4,1}
   &=&  \int_{E^N}  \left( \prod_{i=1}^\ell a_i (X) \right)  \, 
 p_0(dX_1) \, \dots \, p_0(dX_N),
 \end{eqnarray*}
 with 
 \begin{displaymath}
 a_i = a_i(X) := \int_{E} \varphi_i(w) \, S^\infty _t(\hat \mu^N_X)(dw),
 \qquad i=1, \dots, \ell\,,
 \end{displaymath}
and similarly
 \begin{eqnarray*}
 \TT_{4,2} 
 &=& \left\langle \left( S^\infty _t  (p_0) \right)^{\otimes \ell}, 
  \varphi  \right\rangle
 =  \int_{E^N}  \left( \prod_{i=1}^\ell b_i  \right)\, 
 p_0(dX_1) \, \dots \, p_0(dX_N),
 \end{eqnarray*}
 with 
 {
 \begin{displaymath}
 b_i := \int_{E} \varphi_i(w) \, S^\infty _t(p_0)(dw), \qquad i=1, \dots, \ell.  
 \end{displaymath}}
A small calculation gives
 \begin{displaymath}
 \TT_4 \le 
 \sum_{i=1} ^\ell \left( 
 \prod_{j \neq i} \| \varphi_j \|_{L^\infty(E)} \right) \, 
 \int_{E^N} \left| a_i(X) - b_i \right| \, p_0(dX_1) \, \dots \, p_0(dX_N)\,,
 \end{displaymath}
 and finally, using {\sc (H2)}, for any $ 1 \le i \le \ell$ 
 \begin{eqnarray*} \left| a_i(X) - b_i \right| &:=& \left| \int_{E}
     \varphi_i(w) \, \left(S^\infty _t(p_0)(dw) - S^\infty_t(\hat \mu^N
       _X)(dw) \right)\right| \\  
  &\le& C_T ^{\infty,w} \, \| \varphi_i
   \|_{\FF_3} \, \mbox{dist}_{\FF_3 '}\left(p_0, \hat \mu^N _X\right)\,,
 \end{eqnarray*}
  which completes the proof of $\TT_4$, and also the proof of
  Theorem~\ref{thm:main} under the   four hypotheses on the involved
  semigroups. 
 
 \subsection{Applications of the abstract theorem: the Boltzmann equation} 

In this section shall se how the abstract theorem can be applied to
the Boltzmann equation with bounded collision rates. In that case, the
result of Grünbaum~\cite{Grunbaum} can be applied, and so the only new
result that can be deduced from the abstract theorem are the explicit
error bounds.

Some other examples are treated in~\cite{MiMoWe2011arXiv} and
in~\cite{MiMo2010arXiv}, for example
 \begin{itemize}
 \item The McKean-Vlasov equation
 \item The Boltzmann equation with certain classes of  force fields.
 \item The Boltzmann equation for ({\em e.g.}) hard spheres
  \end{itemize}
The last example, which is treated in~\cite{MiMo2010arXiv}, actually
requires some rather technical modifications of the abstract theorem
to handle weighted spaces. All details of this are given
in~\cite{MiMo2010arXiv}.  

The objectiv is then to derive the Boltzmann equation,
\begin{eqnarray*}
\left\{
  \begin{array}{l}
  \partial_t p_t = Q(p_t,p_t)\\
  p_{t=0} = p_0 \,,   
  \end{array}\right.
\end{eqnarray*}
in this case with  with $ p_t\in \PP(\R^3)$, and with the right hand
side defined by 
\begin{eqnarray*}
 \langle Q(p,p), \varphi \rangle := \int_{E^2\times S^2}\underbrace{ \gamma (|w_1 -
 w_2|) \, b(\theta)}_{B(w_1-w_2,\theta)} \, (\phi (w^*_2) - \phi (w_2))
 \, d\sigma \, p(dw_1) \, 
 p(dw_2)  \,.
\end{eqnarray*}
which should hold for any $ \varphi \in C_0(\R^d)$, for any $ p \in
P(\R^d)$, with  
\begin{eqnarray*}
 w_1 ^* = {w_1+ w_2 \over 2} + {|w_2 - w_1|\over 2}\, \sigma, \qquad 
 w_2 ^* = {w_1+ w_2 \over 2} - {|w_2 - w_1|\over 2}\, \sigma.
\end{eqnarray*}
And in order that the collision rate be bounded, we require
$B(w_1-w_2,\theta)$ to be bounded.

The Markov processes on $(\R^3)^N$ are constructed as in the paper by
Grünbaum:
\begin{itemize}
\item[-] For all pairs of indices $ i'\ne j'$ draw $
  T_{i',j'}$ from an exponential distribution with parameter
  $ \gamma(|v_{i'}-v_{j'}|)$
  
  {\em i.e.} $ \mathbb{P}( T_{i',j'} > t) = \exp(-\gamma t)$
\item[-] Let $ T_1 = \min( T_{i',j'} )$ and $ (i,j) =
  (i',j')$
\item[-] Draw $ \sigma\in S^2$ according to law $
  b(\theta_{i,j})$ where $ \cos\theta_{i,j} = \sigma\cdot
  \frac{v_i-v_j}{|v_i'v_j|}$ 
\item[-] The new state after collision at time $T_1$ becomes
  \begin{equation*}
    V^* = V^*_{ij} = R_{ij,\sigma}V = (v_1, ..., v^*_i, ...., v^*_j, ... ,
    v_N),
  \end{equation*}
  with
  \begin{equation*}
    v^*_i = {v_i + v_j \over 2}+ \sigma \, {|v_i - v_j|\over 2}, \quad v^*_j
    = {v_i + v_j\over 2}- \sigma \, {|v_i - v_j|\over 2}
  \end{equation*}
\end{itemize}
The Markov process $ \VV_t$ is constructed by repeating
       the steps above {\em but} with time  rescaled with $ N$ so that
       each coordinate jumps one time  per unit time, on averge.
The law of $ \VV_t$ is denoted  $p^N_t$, and the corresponding
semigroup $ S_t^N$, the dual semigroup $ T_t^N$ and its generator $ G^N$.
 
  The master
 equation on the law $p^N_t$ is given in dual form by
 \begin{displaymath}
   \partial_t \langle p^N_t,\varphi \rangle = \langle p^N_t, G^N
   \varphi \rangle\,,  
 \end{displaymath}
 with 
 \begin{equation*}
   (G^N\varphi) (V) = {1 \over N}\sum_{i,j= 1}^N \gamma(|v_i-v_j|)  
   \int_{\mathbb{S}^{d-1}} b(\theta_{ij}) \, \left[\varphi^*_{ij} -
     \varphi\right] \, d\sigma\,,
 \end{equation*}
 where $ \varphi = \varphi(V)$, $ \varphi^*_{ij}= \varphi(V^*_{ij})$. 
 
 \begin{thm}
\label{thm:boltzmann}
 Assume that $ p_0\in \PP(\R^d)\cap M^1(\R^d;\langle v\rangle^{d+5})$, 
$ p_0^N = \underbrace{p_0\otimes \cdots \otimes
     p_0}_{N\;\mbox{times}}$. Let $ p_t^N = S_t^N(p_0^N)$ be the
   solution of the $N$-particle master equation and
   $p_t=S^{\infty}_t(p_0)$ the solution of the Boltzmann equation.
 Then there is a constant $C(k,\ell)$, depending only on $k$ and
 $\ell$ and $a>0$  such that for $ N>2\ell$, $ 0\le t \le T$, and all 
   \begin{displaymath}
     \varphi = \varphi_1 \otimes \varphi_2 \otimes \dots \otimes \,
   \varphi_\ell \in (C(\R^d) \cap \mbox{{\em Lip}}(\R^d))^{\otimes \ell},
   \end{displaymath}
 we have
   \begin{align*} \label{eq:cvgHS} \sup_{[0,T]}&\left| \left \langle \left(
           S^N_t(p_0 ^N) - \left( S^\infty _t (p_0) \right)^{\otimes N}
         \right), \varphi \right\rangle \right| \\ 
 &\le C(k,\ell) \, \Bigg[
     \frac{\|\varphi\|_{L^\infty} \, \|p_0\|_{M^1}}{N} + 
     e^{a T} \, \frac{\|\varphi\|_{\mbox{{\scriptsize {\em Lip}}}(R^{d\ell})}}{N} +
     \frac{\| \varphi \|_{\mbox{{\scriptsize {\em
               Lip}}}(R^{d\ell})} \, \|
       p_0\|_{M^1_{d+5}}}{N^{1/(d+4)}}\Bigg].
   \end{align*}
  \end{thm}

\noindent
{\sc Proof:} The statement of the theorem is a reformulation of
Theorem~\ref{thm:main}, and the proof is carried out by choosing the
spaces $\FF_1, \FF_2, \FF_3$ and verifying that the hypotheses 
{\sc  (H1)} ... {\sc  (H14)} hold. And this can be done with 
 $ \FF_1=\FF_2 = C_0(\R^d)$ and $ \FF_3= \mbox{{\em Lip}}(\R^d)$.
 
It follows that propagation of chaos holds, at least for this
   kind of initial data. 

\noindent
{\sc Proof of {\sc (H1)}}. We want to show that there exists $ C_1>0$ such that
\begin{displaymath}
  \label{eq:H1BddBoltzmann}
  \forall \, \Phi \in C^{1,1}(M^1), \quad 
  \left\| \left( G^N \, \pi^N - \pi^N \, G^\infty \right)  \,  
    \Phi \right\|_{L^\infty(E^N)} \le {C_1 \over N}   \| \Phi \|_{C^{1,1}(M^1)}.
\end{displaymath} 
For  $ \Phi \in C^{1,1}(M^1)$, set $ \phi = D\Phi[\hat\mu^N_V]$ and compute
\begin{align*}
  G^N &(\Phi \circ \hat \mu^N_V ) 
  = {1 \over 2N} \sum_{i,j= 1}^N \gamma(|v_i-v_j|) 
  \int_{\mathbb{S}^{d-1}} b(\theta_{ij})  \left[ \Phi (
    \hat \mu^N_{V^*_{ij}}) - \Phi ( \hat \mu^N_V)\right] \, d\sigma \\
  &= {1 \over 2N} \sum_{i,j= 1}^N \gamma(|v_i-v_j|)
  \int_{\mathbb{S}^{d-1}} b(\theta_{ij}) \, 
  \langle \hat \mu^N_{V^*_{ij}} -\hat \mu^N_V, \phi \rangle\, d\sigma  
  \qquad\quad (= I_1(V)) \\
  &+ {1 \over 2N} \sum_{i,j= 1}^N \gamma(|v_i-v_j|) 
  \int_{\mathbb{S}^{d-1}} {\mathcal O} \left( 
    \|\Phi \|_{C^{1,1}} \, 
    \left\|\hat \mu^N_{V^*_{ij}} -\hat \mu^N_V \right\|_{M^1}^{2} \right)\, d\sigma
  \quad (= I_2(V)).
\end{align*}
The first term, $I_1$ is estimated with (recall that $\hat{\mu}_V =
\frac{1}{N}\sum_{j=1}^N \delta_{v_j}$)
\begin{align*} 
  I_1 =& {1 \over 2N^2}\sum_{i,j= 1}^N \gamma(|v_i-v_j|) 
  \int_{\mathbb{S}^{d-1}} b(\theta_{ij}) \, 
  \left[\phi(v^*_i) + \phi(v^*_j) -\right.  \\
  &\hspace{0.5\textwidth} \left.\phi(v_i) - \phi(v_j)\right]\, d\sigma \\
  &={1 \over 2}  \int_{\R^d} \int_{\R^d} \int_{\mathbb{S}^{d-1}} \gamma(|v-w|)\, 
  b(\theta) \, 
  \left[\phi(v^*) + \phi(w^*) -  \phi(v) - \phi(w) \right] \, 
  \hat\mu^N_V (dv) \, \hat\mu^N_V (dw) \, d\sigma \\
  &=  \left\langle Q(\hat\mu^N_V,\hat\mu^N_V), \phi \right\rangle 
  = \left(G^\infty \Phi\right) (\hat\mu^N_V)\,,
\end{align*}
and the second, $I_2 (V) $ as
\begin{eqnarray*} 
  I_2 (V) 
  =& {1 \over 2N}\sum_{i,j= 1}^N \gamma(|v_i-v_j|) 
  \int_{\mathbb{S}^{d-1}} {\mathcal O} 
  \left( \|\Phi \|_{C^{1,1}} \, 
    \left( {4 \over N} \right)^{2} \right)\, d\sigma \\
  &\le C \, \| \gamma \|_\infty \, {\|\Phi \|_{C^{1,1}} \over N}  \, 
  \left( \sum_{i,j = 1}^N {1 \over N^2} \right)
  \le C \, {\|\Phi \|_{C^{1,1}}  \over N}.
\end{eqnarray*}
And together these yield the desired estimate.

\medskip
\noindent
{\sc Proof of {\sc (H2)}}. 
 Set $ k=\theta=1$. We will show that for all  $ \mu, \, \mu' \in P(\R^d)$
 and for any $ T>0$, there exists $ C_T>0$ such that
\begin{equation}
\label{eq:h2boltz}
  \sup_{t \in [0,T]} \Big\| S^\infty _t(\mu') - S^\infty_t(\mu)  
  -  {\mathcal LS}_t^\infty[\mu] (\mu' - \mu) \Big\|_{M^1} 
  \le C_T\, \| \mu' - \mu \|_{M^1}^2, 
   \end{equation}
    where $ {\mathcal LS}^\infty_t[\mu]$ is the linear semigroup associated
   the solution $ S^\infty_t \mu$. Hence  
   {\sc (H2)} holds with $ \FF_2 = \FF_1 = C_0(\R^d)$,
   $\FF_2' = M^1 (\R^d)$  when $ k=\theta=1$.
  To prove~(\ref{eq:h2boltz}) consider
  \begin{eqnarray*}
 \partial_t f_t &=& Q(f_t,f_t), \qquad f_0 = \mu, \vspace{0.3cm} \\ \displaystyle
 \partial_t g_t &=& Q(g_t,g_t), \qquad g_0 = \mu', \vspace{0.3cm} \\
 \displaystyle 
 \partial_t h_t &=& \tilde Q(f_t,h_t) := Q(f_t,h_t) + Q(h_t,f_t), 
 \qquad h_0 = g_0 - f_0 = \mu' - \mu.    
  \end{eqnarray*}
The solutions to these equations, $f$, $g$, and $h$, and $\phi \equiv
f-g-h$ is the remainder term in~(\ref{eq:h2boltz}), and this can be
estimated by a Gronwall argument to give the estimate, with $C_T\sim
e^{C T}$. 

\medskip
\noindent
{\sc Proof of {\sc H3}}:  Take$ \FF_3 = \mbox{Lip}(\R^d)$. The
Wasserstein (or Tanaka) distance between two measures is defined as
\begin{equation*}
  W_1(\mu,\nu) = \inf\limits_{\gamma\in\Gamma(\mu,\nu)} \int_{E\times
    E} |x-y|\,d\gamma(x,y)
    =\inf \Ee\left[ |X-Y|\right]\,.
\end{equation*}
  Here $ \Gamma(\mu,\nu)$ is the collection of $
  \gamma\in\PP(E\times E)$ 
  such that 
  $ \mu = \int_{E} \gamma(\cdot,dy)$ and
  $ \nu = \int_{E} \gamma(dx,\cdot)$.
 \medskip
An equivalent definition is 
\begin{equation*}
 W_1(\mu,\nu) = \sup \left\{
 \left|
 \int \phi(x) ( \mu(dx)-\nu(dx))
 \right|
 \;:\; \| \phi \|_{\mbox{\em Lip}} \le 1
 \right\}  \,.
\end{equation*}
Tanaka~\cite{T1,T5} proved that if $ p_t^1, p_t^2$ are the solution so of
the Boltzmann equation with initial data 
 $ p_0^1, \, p_0^2 \in P(\R^d)$, then 
 \begin{equation*}
 \sup_{[0,T]} W_1(p^1_t,p^2_t) \le W_1(p^1_0,p^2_0),\,.  
 \end{equation*}
Now  {\sc (H3)}
   \begin{displaymath}
     \forall  \, \mu, \nu \in P(E)  \qquad 
     \sup_{[0,T]} \mbox{dist}_{\FF'_3} 
     \left( S^\infty _t(\mu), S^\infty _t(\nu) \right) \le
     C_T^{\infty,w} \, \mbox{dist}_{\FF'_3} (\mu,\nu) \,,
   \end{displaymath}
 with $\FF_3= \mbox{Lip}(\R^d)$, follows immediatley from Tanaka's result.

\medskip
\noindent
{\sc Proof of {\sc H4}}: Here we need to check that the
   dual of $ \FF_3 = \mbox{Lip}(\R^d)$,  $ \FF_3 '$ satisfies
\begin{displaymath}
     \left\| \pi^N (\Phi) \right\|_{\FF_3 ^N} =
   \left\| \Phi \circ \hat \mu^N _X \right\|_{\FF_3 ^N} 
   \le C_\pi \, \left\| \Phi \right\|_{C^{0,1}(\FF_3 ')}.
 \end{displaymath}
 For any $ \Phi \in C^{0,1}(\FF_3')$,
 \begin{eqnarray*}
 \left\| \pi^N \Phi \right\|_{\mbox{Lip}((\R^d)^N)} &\le& 
 \sup_{X\neq Y \in (\R^d)^N} \frac{\left|\Phi(\hat \mu^N _X) - \Phi(\hat \mu^N
   _Y)\right|}{|X-Y|} \\ 
 &\le& \| \Phi\|_{C^{0,1}(\FF_3')} \, 
 \frac{W_1(\hat \mu^N _X,\hat \mu^N _Y)}{|X-Y|} \le 
 C \, \| \Phi\|_{C^{0,1}(\FF_3')}.
 \end{eqnarray*}
This implies that {\sc (H4)}, and hence all four hypotheses are
satisfied, and Theorem~\ref{thm:boltzmann} is a consequence of
Theorem~\ref{thm:main}. 

 \subsection{Other examples and comments}
 
Another model that is covered by Theorem~\ref{thm:main} is the
McKean-Vlasov system~\cite{McK3}. Here the $N$-particle system is
defined as
 \begin{displaymath}
 d x^i_t = \sigma_i \, dB^i_t + F^N(x_i, \hat\mu^{N-1}_{\hat X^i})\, dt\,, \qquad 1 \le i \le N\,,
 \end{displaymath}
 with $ \hat X^i:= ( x^1, ...,x^{i-1},x^{i+1}, ... , x^N)$ and
 $ F^N : \R^d
 \times P(\R^d) \to \R^d$.
 The nonlinear McKean-Vlasov
 equation on $ P(\R^d)$ defined by
 \begin{displaymath}{\partial p \over \partial t} = Q(p_t),
   \quad p(0) = p_0 \quad\hbox{in}\quad P(\R^d),
 \end{displaymath}
 with
 \begin{equation*}
 Q(\rho) = {1 \over 2}\sum_{\alpha,\beta =1}^d \partial^2_{\alpha \beta}
 (\Sigma_{\alpha\beta}\, \rho) - \sum_{\alpha=1}^d \partial_\alpha
 (F_\alpha(x,\rho) \, \rho).   
 \end{equation*}
In this case, the hypotheses {\sc (H1)} to {\sc (H4)} can be verified with 
$ \FF_1=H^2, \FF_2=H^{s+2}, \FF_3 = Lip(\R^d)$

But the abstract theorem presented here does not cover e.g the
Boltzmann equation for hard spheres, {\em i.e.} the case that Grünbaum
attempted to solve. A more detailed analysis, involving weighted
spaces, is required for that. A proof is given
in~\cite{MiMo2010arXiv}.

Another important result in~\cite{MiMo2010arXiv} is that in some cases
all estimates can be carried out uniformly in time (contrary to the
estimate above, which involves constants that grow exponentially with
the time interval). It is not at all obvious that such a result could
be true, considering the calculations carried out in
Section~\ref{sec:Kac}. For large times the exponential $e^{tL}$ will
be dominated by large powers of $L$, and for any fixed $N$, the same
variables must be reused many times, potentially creating correlations
that remain also when $N$ increases. For the Boltzmann equation and
the related $N$-particle systems, the stationary measures to the
$N$-particle systems are themselves chaotic, and this may help getting
the uniform estimates.

However, the model of flocking described in Section~\ref{sec:BDG} does not
have this property. It is a ``pair interaction driven master
equation'' which are  defined in~\cite{CarlenDegondWennberg2011},
where it is also proven that propagation of chaos holds for all
times, but that the stationary states for the $N$-particle systems are
{\em not} chaotic. Another model studied
in~\cite{CarlenDegondWennberg2011} is called a ``choose the leader
model''. In that model a pair interacts in such a way that one of the
two particles (randomly chosen in the pair) tries to take the other
particle's velocity, but makes a random error. That is also a pair
interaction driven master equation, and in this case some calculations
can be carried out rexplicitly, and in particular one can find
explicit expressions for the marginal distributions. These expressions
show that the stationary states are not chaotic.

Propagation of chaos is an important concept, and many questions
remain open, most notably the question of propagation of chaos for a
deterministic particle system and a rigorous derivation of the
Boltzmann equation, valid over a macroscopic time interval. I hope
that these notes have given some flavour of this and recommend the
reader to look in the litterature for many more results. Some relevant
references are~\cite{CarlenCarvalhoLeRouxLossVillani2010,AmmariNierIHP2008,
GrahamMeleardESIAM1999,JourdainMeleardIHP1998,CaprinoPulvirentiWagnerSIAM1998,
SznitmanSaintFlur1989,CaprinoDeMasiPresuttiPulvirentiJSP1989,SznitmanZWVG1984,
McK3}.

\subsection*{Acknowledgment}

I would like to express my gratitude to the orgnization committee of
5th Summer School on  ''METHODS AND MODELS OF KINETIC THEORY'' for
giving me the opportunity to give this series of lectures. I would
also like to thank my co-authors in the papers that form a bases for
the notes: Eric Carlen, Pierre Degond, Johan Henriksson, Torbjörn Lundh, Stéphane Mischler, Clément Mouhot.


\bibliographystyle{acm}

\end{document}